\documentclass[12pt,a4paper]{article}
\usepackage[dvipdfmx]{graphicx}
\usepackage{latexsym}
\usepackage{amssymb}
\usepackage{amsmath}
\usepackage{amscd}
\usepackage{amsthm}
\usepackage{amsfonts}
\usepackage{mathrsfs}
\usepackage{enumerate}
\usepackage{bm}
\usepackage[usenames]{color}
 \makeatletter
    
    \@addtoreset{equation}{section}
  \makeatother
\newtheorem{Th}{Theorem}[section]
\newtheorem{Co}[Th]{Corollary}
\newtheorem{Lem}[Th]{Lemma}
\newtheorem{Pro}[Th]{Proposition}
\theoremstyle{definition}
\newtheorem{Exa}[Th]{Example}
\newcommand{\demo}{\par\noindent{\it Proof. \/}\ }
\newcommand{\enD}{\hfill $\Box$\vspace{3truemm} \par}

\newcommand{\R}{{\mathbb R}}
\newcommand{\lon}{\longrightarrow}

\begin{document}
\title{Primitivoids and inversions of plane curves
 }
\author{Shyuichi IZUMIYA and Nobuko TAKEUCHI }

\date{\today}
\maketitle
\begin{abstract}
The pedal of a curve in the Euclidean plane is a classical subject
which has a singular point at the inflection point of the original curve.
The primitive of a curve is a curve given by the inverse construction for
making the pedal.
We consider relatives of the primitive of a plane curve
which we call {\it primitivoids}.
We investigate the relationship of primitivoids and pedals of plane curves.
\end{abstract}
\renewcommand{\thefootnote}{\fnsymbol{footnote}}

\section{Introduction} 
For a regular curve in the Euclidean plane, the {\it pedal} of a plane curve is the locus of the points on the tangent line of the curve, which are given by the projection image of the position vector of the curve
along the normal direction.
A point on the pedal determines the tangent line of the original curve
at the corresponding point.
Therefore, the passage to the pedal is sometimes called the formation of a \lq\lq derivative\rq\rq.
The inverse operation is called the formation of a\lq\lq primitive\rq\rq.
The {\it primitive} of a curve in the plane is the envelope of the normal lines to its 
position vectors at their ends (cf.\cite[pp. 91]{Arnold}).
It is known that  the singularities of the pedal correspond to 
the inflection points of the original curve.
In this paper we introduced the notion of the anti-pedal of a curve
whose singularities also corresponding to the inflection points of the original curve.
Moreover, we show that the primitive is equal to the anti-pedal of the inversion image
of the original curve. 
In \cite{I-T} we defined the \textit{$\phi$-pedaloid} of a curve as
the locus of the points on the line with a constant angle $\phi$ to the tangent line of the original curve
which are given by the projection image of the position vector of the original curve along the normal
direction of the line.
Therefore, the $\pi/2$-pedaloid of a curve is the pedal of the curve.
Then the family of $\phi$-pedaloid is called the \textit{pedaloids}.
The notion of pedaloids is an analogous notion of evolutoids in \cite{G-W}.
In \cite{I-T} we investigated the relation between pedaloids and evolutoids of a curve.
In this paper we also introduce the notion of primitivoids of a plane curve, which are relatives of
the primitive.
There are two ways to define the primitivoids.
One is called a {\it parallel primitivoid}.
For $r\in \R\setminus \{0\},$ the {\it $r$-parallel primitivoid} of a curve is defined
to be the envelope of the family of normal lines to its position vectors at their $r$ times the ends.
Another is called a {\it slant primitivoid}. 
For $\phi\in \R,$ the {\it $\phi$-slant primitivoid} of a curve is
defined to be the envelope of the family of lines
with the constant angle $\phi$ to the position vectors at their ends.
We give the exact parametrizations of both of the $r$-parallel primitivoid and the $\phi$-slant primitivoid of the curve.
In \cite{I-T} we have drawn the pictures of pedaloids which have complicated shapes.
By definition, the shape of the primitivoids might be also complicated.
However, we show that the shapes of all primitivoids of a given curve are similar (cf. Theorems 3.1 and 4.4).
We remark that the primitivoids and the evolutoids of a circle around the origin are
the same (i.e. concentric circles).
However, these are quite different for general curves.
The notion of primitivoids is deeply related to the notion of
apertures of plane curves (cf. \cite{K-N}).
Finally we extend the notion of primitivoids for a kind of singular curves (i.e. frontals) in
the plane.
\section{Pedals, anti-pedals and primitives}
In this section we quickly review the properties of pedals and primitives.
Moreover, we introduce the notion of anti-pedals which  plays an important role in this paper.
Let $\bm{\gamma}:I\lon \R^2$ be a unit speed plane curve, where $\R^2$ is the Euclidean plane with the
canonical scaler product $\langle \bm{a},\bm{b}\rangle=a_1b_1+a_2b_2$ for $\bm{a}=(a_1,a_2),\bm{b}=(b_1,b_2)\in \R^2.$
The {\it norm} of $\bm{a}$ is defined to be $\|\bm{a}\|=\sqrt{\langle \bm{a},\bm{a}\rangle}=\sqrt{a_1^2+a_2^2}.$
Then we have the Frenet formulae:
\[
\left\{
\begin{array}{ll}
\bm{t}'(s)=\kappa(s)\bm{n}(s), \\
\bm{n}'(s)=-\kappa (s)\bm{t}(s),
\end{array}
\right.
\]
where $\bm{t}(s)=\bm{\gamma}'(s)$ is the {\it unit tangent vector}, $\bm{n}(s)=J\bm{t}(s)$ is the {\it unit normal vector} and
$\kappa (s)=x_1'(s)x_2''(s)-x_1''(s)x_2'(s)$ is the {\it curvature} of $\bm{\gamma}(s)$, where  $J= \begin{pmatrix} 0 & -1 \\ 1 & 0 \end{pmatrix}$ and $\bm{\gamma}(s)=(x_1(s),x_2(s)).$
The {\it pedal curve} of $\bm{\gamma}$ is defined to be
${\rm Pe}_{\bm{\gamma}}(s)=\langle \bm{\gamma}(s),\bm{n}(s)\rangle \bm{n}(s)$ (cf. \cite[Page 36]{Bru-Gib}).
Since
\[
{\rm Pe}_{\bm{\gamma}}'(s)=-\kappa (s)(\langle \bm{\gamma} (s),\bm{t}(s)\rangle\bm{n}(s)+\langle \bm{\gamma}(s),\bm{n}(s)\rangle\bm{t}(s)),
\]
the singular points of the pedal of $\bm{\gamma}$ is the point $s_0$ where  $\bm{\gamma}(s_0)=\bm{0}$ or
$\kappa (s_0)=0$ (i.e. the inflection point).
If we assume that $\bm{\gamma}$ does not pass through the origin, the singular points of the pedal 
${\rm Pe}_{\bm{\gamma}}$ are the inflection points of $\bm{\gamma}.$
Therefore we assume that $\bm{\gamma}$ does not pass through the origin.
By definition, ${\rm Pe}_{\bm{\gamma}}(s)$ is the point on the tangent line through $\bm{\gamma}(s)$, which is given by the projection image of $\bm{\gamma}(s)$ of the normal direction. 
Thus, ${\rm Pe}_{\bm{\gamma}}(s)-\bm{\gamma}(s)$ generates the tangent line at $\bm{\gamma}(s).$
The {\it primitive} of a curve $\bm{\gamma}$ in the plane is the envelope of the normal lines to its 
position vectors $\bm{\gamma}(s)$ at their ends (cf.\cite[pp. 91]{Arnold}).
 For a unit speed plane curve $\bm{\gamma}: I\lon \R^2\setminus\{\bm{0}\},$ we define a family of functions
 $H:I\times (\R^2\setminus \{\bm{0}\})\lon \R$ by $H(s,\bm{x})=\langle \bm{x}-\bm{\gamma}(s),\bm{\gamma}(s)\rangle.$
 For any fixed $s\in I,$ $h_s(\bm{x})=H(s,\bm{x})=0$ is an equation of the line through $\bm{\gamma}(s)$ and orthogonal to the position vector $\bm{\gamma}(s).$
The envelope of the family of the lines $\{h_s^{-1}(0)\}_{s\in I}$ is the primitive of $\bm{\gamma}.$
Since $\partial H/\partial s(s,\bm{x})=\langle \bm{x}-2\bm{\gamma}(s),\bm{t}(s)\rangle,$ the primitive ${\rm Pr}_{\bm{\gamma}}:I\lon \R^2\setminus \{\bm{0}\}$ of 
$\bm{\gamma}$ is given by
\[
{\rm Pr}_{\bm{\gamma}}(s)=2\bm{\gamma}(s)-\frac{\|\bm{\gamma}(s)\|^2}{\langle\bm{n}(s),\bm{\gamma}(s)\rangle}\bm{n}(s).
\]
The pedal and the primitive of a regular curve generally have singularities.
Then the primitive of the pedal and the pedal of the primitive for the regular part is well defined, respectively. Suppose that ${\rm Pe}_{\bm{\gamma}}$ and ${\rm Pr}_{\bm{\gamma}}$ are regular curves.
Then we have ${\rm Pr}_{{\rm Pe}_{\bm{\gamma}}}(s)={\rm Pe}_{{\rm Pr}_{\bm{\gamma}}}(s)=\bm{\gamma}(s).$
These arguments for singular cases can be naturally interpreted by using the notion of frontals in \S 3.
\par
On the other hand, it is known that the pedal is given as the envelope of a family of circles as follows:
Let $G:I\times \R^2\lon \R$ be a function defined by
\[
G(s,\bm{x})=\left\|\bm{x}-\frac{1}{2}\bm{\gamma}(s)\right\|^2-\frac{1}{4}\|\bm{\gamma}(s)\|^2=\langle \bm{x},\bm{x}-\bm{\gamma}(s)\rangle.
\]
If we fix $s_0\in I,$ $G(s_0,\bm{x})=0$ is the equation of the circle with the center $\frac{1}{2}\bm{\gamma}(s_0)$ which passes through the origin.
\par
For a fixed $s\in I,$ $g_s(\bm{x})=G(s,\bm{x})=0$ is an equation of a circle through the origin.
Therefore the inversion image of it is a line.
We define the {\it inversion} $\Psi :\R^2\setminus \{\bm{0}\}\lon \R^2\setminus \{\bm{0}\}$ at the origin with respect to the unit circle by
$\displaystyle{\Psi (\bm{x})=\frac{\bm{x}}{\|\bm{x}\|^2}}.$ Then we have $\Psi({g_s^{-1}(0)})=\{\bm{x}\ |\ \langle \bm{x},\bm{\gamma}(s)\rangle =1\}.$
Thus, we define a family of functions
$F:I\times (\R^2\setminus\{\bm{0}\})\lon \R$ by
$F(s,\bm{x})=\langle \bm{x},\bm{\gamma}(s)\rangle -1.$
Then we have
\[
\frac{\partial F}{\partial s}(s,\bm{x})=\langle \bm{x},\bm{t}(s)\rangle,
\]
so that $F=\partial F/\partial s=0$ if and only if 
\[
\bm{x}=\frac{1}{\langle \bm{\gamma}(s),\bm{n}(s)\rangle}\bm{n}(s).
\] 
We define a mapping ${\rm APe}_{\bm{\gamma}}:I\lon \R^2\setminus \{\bm{0}\}$ by
\[
{\rm APe}_{\bm{\gamma}}(s)=\frac{1}{\langle \bm{\gamma}(s),\bm{n}(s)\rangle}\bm{n}(s),
\]
which is called an {\it anti-pedal curve} of $\bm{\gamma}.$
By definition or a straightforward calculation, we have 
$\Psi \circ {\rm APe}_{\bm{\gamma}}={\rm Pe}_{\bm{\gamma}}$ and $\Psi \circ {\rm Pe}_{\bm{\gamma}}={\rm APe}_{\bm{\gamma}}.$
We have the following proposition.
\begin{Pro} For any unit speed plane curve $\bm{\gamma}:I\lon \R^2\setminus\{\bm{0}\},$
we have
\[
{\rm Pr}_{\bm{\gamma}}(s)={\rm APe}_{\Psi\circ\bm{\gamma}}(s)\ \mbox{and}\ {\rm Pr}_{\Psi\circ\bm{\gamma}}(s)={\rm APe}_{\bm{\gamma}}(s).
\]
\end{Pro}
\demo
We do not use the parametrization of the primitive and use the properties of the envelope of the family of lines $\{h_s^{-1}(0)\}_{s\in I}$.
Since $H(s,\bm{x})=\langle \bm{x},\bm{\gamma}(s)\rangle-\|\bm{\gamma}(s)\|^2,$
$H(s,\bm{x})=0$ if and only if $\langle \bm{x},\Psi\circ \bm{\gamma}(s)\rangle =1.$
Therefore, the envelope of the family of lines $\{h_s^{-1}(0)\}_{s\in I}$ is equal to the anti-pedal of $\Psi\circ \bm{\gamma}.$
This means that ${\rm Pr}_{\bm{\gamma}}(s)={\rm APe}_{\Psi\circ\bm{\gamma}}(s).$ 
Since $\Psi\circ \Psi =1_{\R^{2}\setminus \{\bm{0}\}},$ we have
\[
{\rm APe}_{\bm{\gamma}}(s)={\rm APe}_{\Psi\circ\Psi\circ\bm{\gamma}}(s)={\rm Pr}_{\Psi\circ\bm{\gamma}}(s).
\]
This completes the proof.
\enD
\par
Since $\Psi$ is a diffeomorphism, the pedal ${\rm Pe}_{\bm{\gamma}}$ and the anti-pedal
${\rm APe}_{\bm{\gamma}}=\Psi\circ {\rm Pe}_{\bm{\gamma}}$ have the same singularities,
which correspond to the inflection points of $\bm{\gamma}.$
Therefore, the singularities of the primitive ${\rm Pr}_{\bm{\gamma}}={\rm APe}_{\Psi\circ\bm{\gamma}}$
correspond to the inflections of the inversion curve $\Psi\circ\bm{\gamma}.$
Therefore, we now calculate the curvature of $\Psi\circ\bm{\gamma}.$
\begin{Pro}
For a unit speed curve $\bm{\gamma}:I\lon \R^2\setminus\{\bm{0}\},$
the curvature of $\Psi\circ \bm{\gamma}$ at $s\in I$ is
$\kappa_{\Psi\circ\bm{\gamma}}(s)=-\kappa (s)\|\bm{\gamma}(s)\|^2-2\langle \bm{\gamma}(s),\bm{n}(s)\rangle.$
\end{Pro}
\demo
We denote that $\widetilde{\bm{\gamma}}(s)=\Psi\circ \bm{\gamma}(s).$
Since $\widetilde{\bm{\gamma}}'=\frac{\|\bm{\gamma}\|^2\bm{t}-2\langle\bm{\gamma},\bm{t}\rangle\bm{\gamma}}{\|\bm{\gamma}\|^4},$ $\|\widetilde{\bm{\gamma}}'\|=1/\|\bm{\gamma}\|^2.$
Let $\sigma$ be the arc-length parameter of $\widetilde{\bm{\gamma}}.$
Then $d\sigma/ds=1/\|\bm{\gamma}\|^2.$
We also have
\[
\widetilde{\bm{t}}=\frac{\widetilde{d\bm{\gamma}}}{d\sigma}=\frac{\|\bm{\gamma}\|^2\bm{t}-2\langle \bm{\gamma},\bm{t}\rangle\bm{\gamma}}{\|\bm{\gamma}\|^2}=
\frac{(\langle\bm{\gamma},\bm{n}\rangle^2-\langle\bm{\gamma},\bm{t}\rangle^2)\bm{t}-2\langle\bm{\gamma},\bm{t}\rangle\langle\bm{\gamma},\bm{n}\rangle\bm{n}}{\|\bm{\gamma}\|^2}.
\]
It follows that
\[
\widetilde{\bm{n}}=J\widetilde{\bm{t}}=\frac{2\langle\bm{\gamma},\bm{t}\rangle\langle\bm{\gamma},\bm{n}\rangle\bm{t}+(\langle\bm{\gamma},\bm{n}\rangle^2-\langle\bm{\gamma},\bm{t}\rangle^2)\bm{n}}{\|\bm{\gamma}\|^2}.
\]
By the Frenet formulae and the relation $\bm{\gamma}=\langle\bm{\gamma},\bm{t}\rangle\bm{t}+\langle\bm{\gamma},\bm{n}\rangle\bm{n},$ we have
\[
\kappa_{\Psi\circ\bm{\gamma}}\widetilde{\bm{n}}=\frac{d\widetilde{\bm{t}}}{d\sigma}
=\|\bm{\gamma}\|^2\frac{d\widetilde{\bm{t}}}{ds}
=-(\kappa \|\bm{\gamma}\|^2+2\langle \bm{\gamma},\bm{n}\rangle)\left(\frac{2\langle\bm{\gamma},\bm{t}\rangle\langle\bm{\gamma},\bm{n}\rangle\bm{t}+(\langle\bm{\gamma},\bm{n}\rangle^2-\langle\bm{\gamma},\bm{t}\rangle^2)\bm{n}}{\|\bm{\gamma}\|^2}\right).
\]
This means that $\kappa_{\Psi\circ\bm{\gamma}}(s)=-\kappa (s)\|\bm{\gamma}(s)\|^2-2\langle \bm{\gamma}(s),\bm{n}(s)\rangle.$
\enD
Therefore, we have
$
\kappa_{\Psi\circ\bm{\gamma}}'(s)=-\kappa'(s)\|\bm{\gamma}(s)\|^2,
$
so that the following corollary holds.
\begin{Co} 
For a unit speed curve $\bm{\gamma}:I\lon \R^2\setminus \{\bm{0}\},$
$\kappa '(s_0)=0$ if and only if $\kappa_{\Psi\circ\bm{\gamma}}'(s_0)=0.$
\end{Co}
\noindent
The point $s_0\in I$ with $\kappa'(s_0)=0$ is called a {\it vertex} of $\bm{\gamma}.$
Thus the above corollary asserts that the vertices of $\bm{\gamma}$ and $\Psi\circ \bm{\gamma}$ are the same.
It is known that the vertex of a curve is the point where the osculating circle has
four point contact with the curve.
Since the inversion $\Psi$ is a conformal diffeomorphism, it sends the osculating circle of 
$\bm{\gamma}$ to
the osculating circle of $\Psi\circ \bm{\gamma}.$
Therefore the above result is reasonable. By definition, a point $s_0\in I$ is 
an inflection point of $\bm{\gamma}$ if and only if $\kappa (s_0)=0.$
We also have the following corollary.
\begin{Co}
For a unit speed curve $\bm{\gamma}:I\lon \R^2\setminus \{\bm{0}\},$
$s_0\in I$ is an inflection point of $\Psi\circ\bm{\gamma}$ if and only if
$\kappa (s_0)\|\bm{\gamma}(s_0)\|^2+2\langle \bm{\gamma}(s_0),\bm{n}(s_0)\rangle=0.$
\end{Co}
\noindent
We now consider the geometric meaning of the
condition $\kappa (s_0)\|\bm{\gamma}(s_0)\|^2+2\langle \bm{\gamma}(s_0),\bm{n}(s_0)\rangle=0.$
We assume that $\langle \bm{\gamma}(s_0),\bm{n}(s_0)\rangle\not=0,$
so that $\kappa (s_0)\not=0.$
The center of the curvature at $s_0$ is given by
\[
\bm{x}_0=\bm{\gamma}(s_0)+\frac{1}{\kappa (s_0)}\bm{n}(s_0).
\]
The osculating circle of $\bm{\gamma}$ at $s_0$ is
defined by the equation
\[
\|\bm{x}-\bm{x}_0\|^2=\frac{1}{\kappa (s_0)^2}.
\]
Therefore, the osculating circle contains the origin if and only if
\[
\left\|\bm{\gamma}(s_0)+\frac{1}{\kappa (s_0)}\bm{n}(s_0)\right\|^2=\|\bm{x}_0\|^2=\frac{1}{\kappa (s_0)^2}.
\]
The last condition is equivalent to $\kappa (s_0)\|\bm{\gamma}(s_0)\|^2+2\langle \bm{\gamma}(s_0),\bm{n}(s_0)\rangle=0.$
Therefore, we have shown the following lemma.
\begin{Lem} With the same assumptions as above, an inflection point of $\Psi\circ\bm{\gamma}$ corresponds 
to the point of $\bm{\gamma}$ such that the osculating circle of the point passes through the
origin.
\end{Lem}
We remember that the osculating circle of a curve at a point is the circle which has three points contact 
with the curve at the point.
If the osculating circle passes through the origin, the inversion image is a tangent line
of the inversion curve which has three points contact.
This means that the point is an inflection point of the inversion curve.
Therefore the above result is also reasonable.
\par
We now consider the singularities of primitives.
By Proposition 2.1, we have ${\rm Pr}_{\bm{\gamma}}(s)={\rm APe}_{\Psi\circ \bm{\gamma}}(s)
=\Psi\circ {\rm Pe}_{\Psi\circ \bm{\gamma}}(s).$
Since $\Psi$ is a conformal diffeomorphism,
the singularities of the primitive ${\rm Pr}_{\bm{\gamma}}$ are equal to the singularities
of the pedal ${\rm Pe}_{\Psi\circ \bm{\gamma}}$ of the inversion curve $\Psi\circ \bm{\gamma}.$
Here, $s_0\in I$ is a singular point of ${\rm Pe}_{\Psi\circ \bm{\gamma}}$ if and only if
$\kappa_{\Psi\circ\bm{\gamma}}(s_0)=0.$
Therefore, we have shown the following proposition.
\begin{Pro}
Let $\bm{\gamma}:I\lon \R^2\setminus \{\bm{0}\}$ be a unit speed curve such that
$s_0\in I$ satisfies $\langle \bm{\gamma}(s_0),\bm{n}(s_0)\rangle\not=0.$ 
Then $s_0$ is a singular point of the primitive ${\rm Pr}_{\bm{\gamma}}$ if and only if
$s_0$ is an inflection point of the inversion curve $\Psi\circ \bm{\gamma}.$
Moreover, such a point satisfies the condition that the osculating circle at the point
$\bm{\gamma}(s_0)$ passes through the origin.
\end{Pro}
As an application of the unfolding theory of functions of one-variable, we have the following proposition (cf. \cite[Page 166]{Bru-Gib}).
\begin{Pro}
The pedal $\textrm{Pe}_{\bm{\gamma}}$ of a curve $\bm{\gamma}$ around $\textrm{Pe}_{\bm{\gamma}}(s_0)$ is locally diffeomorphic to
the \textit{ordinary cusp} $C=\{\ (x,y)\in \R\ \mid \ x=t^2,y=t^3\, \}$ if and only if $\kappa (s_0)=0$ and $\kappa '(s_0)\not=0.$ 
\end{Pro}
We say that $\bm{\gamma}(s_0)$ is an {\it ordinary inflection} of $\bm{\gamma}$ if $\kappa (s_0)=0$ and $\kappa '(s_0)\not=0.$ Then we have the following theorem.
\begin{Th} Let $\bm{\gamma}:I\lon \R^2\setminus \{\bm{0}\}$ be a unit speed curve such that
$s_0\in I$ satisfies $\langle \bm{\gamma}(s_0),\bm{n}(s_0)\rangle\not=0.$ 
Then the following conditions are equivalent\/{\rm :}
\par\noindent
{\rm (1)} The primitive ${\rm Pr}_{\bm{\gamma}}$ around ${\rm Pr}_{\bm{\gamma}}(s_0)$ is locally 
diffeomorphic to the ordinary cusp,
\par\noindent
{\rm (2)} the pedal ${\rm Pe}_{\Psi\circ \bm{\gamma}}$ of $\Psi\circ\bm{\gamma}$ around
${\rm Pe}_{\Psi\circ \bm{\gamma}}(s_0)$ is locally diffeomorphic to the ordinary cusp,
\par\noindent
{\rm (3)}  $\kappa (s_0)\|\bm{\gamma}(s_0)\|^2+2\langle \bm{\gamma}(s_0),\bm{n}(s_0)\rangle=0$ and $\kappa'(s_0)\not= 0,$
\par\noindent
{\rm (4)} the point $\bm{\gamma}(s_0)$ is not a vertex $\bm{\gamma}$ and the osculating circle of $\bm{\gamma}$ at $\bm{\gamma}(s_0)$ passes through the origin,
\par\noindent
{\rm (5)} the point $\Psi\circ \bm{\gamma}(s_0)$ is an ordinary inflection of $\Psi\circ\bm{\gamma}.$
\end{Th}
\demo
By Proposition 2.7, conditions (2) and (5) are equivalent.
By Corollaries 2.3 and 2.4, conditions (3) and (5) are equivalent.
By Corollary 2.3 and Lemma 2.5, conditions (4) and (5) are equivalent.
By Proposition 2.1, ${\rm APe}_{\Psi\circ \bm{\gamma}}={\rm Pr}_{\bm{\gamma}}.$
Since $\Psi$ is a diffeomorphism, the pedal ${\rm Pe}_{\Psi\circ \bm{\gamma}}$ around
${\rm Pe}_{\Psi\circ \bm{\gamma}}(s_0)$ is locally 
diffeomorphic to the ordinary cusp if and only if ${\rm APe}_{\Psi\circ \bm{\gamma}}=\Psi\circ {\rm Pe}_{\Psi\circ \bm{\gamma}}$ around
${\rm APe}_{\Psi\circ \bm{\gamma}}(s_0)$ is locally diffeomorphic to the ordinary cusp.
Thus conditions (1) and (2) are equivalent. This completes the proof.
\enD

\section{Parallel primitivoids}
In this section we give a parametrization of the $r$-parallel primitivoid of $\bm{\gamma}$
for $r\in \R.$
We consider a family of functions $H:I\times (\R\setminus \{\bm{0}\})\lon \R$
defined by $H(s,\bm{x})=\langle \bm{x}-r\bm{\gamma}(s),\bm{\gamma}(s)\rangle.$
Then $\{h_s^{-1}(0)\}_{s\in I}$ is a family of lines which are orthogonal to
$\bm{\gamma}(s)$ through $r\bm{\gamma}(s),$ where $h_s(\bm{x})=H(s,\bm{x}).$
By definition, the envelope of this family of lines is the $r$-parallel primitivoid of $\bm{\gamma}.$
Since $\{\bm{t}(s),\bm{n}(s)\}$ is an orthonormal frame along $\bm{\gamma},$
there exist $\lambda, \mu \in \R$ such that
$\bm{x}-r\bm{\gamma}(s)=\lambda \bm{t}(s)+\mu \bm{n}(s).$
We have
\[
\frac{\partial H}{\partial s}(s,\bm{x})=\langle -r\bm{t}(s),\bm{\gamma}(s)\rangle +\langle \bm{x}-r\bm{\gamma}(s), \bm{t}(s)\rangle=\langle \bm{x}-2r\bm{\gamma}(s),\bm{t}(s)\rangle.
\] 
We also have 
\[
\langle \bm{x}-2r\bm{\gamma}(s),\bm{t}(s)\rangle=\langle -r\bm{\gamma}(s)+\lambda \bm{t}(s)+\mu \bm{n}(s),
\bm{t}(s)\rangle=-r\langle \bm{\gamma}(s),\bm{t}(s)\rangle+\lambda,
\]
so that
$\partial H/\partial s(s,\bm{x})=0$ if and only if $\lambda = r\langle \bm{\gamma}(s),\bm{t}(s)\rangle.$
Moreover, $F(s,\bm{x})=0$ if and only if
$\lambda \langle \bm{t}(s),\bm{\gamma}(s)\rangle +\mu\langle \bm{n}(s),\bm{\gamma}(s)\rangle =0.$
Thus $F(s,\bm{x})=\partial H/\partial s(s,\bm{x})=0$ if and only if
\[
0=\langle r\langle \bm{t}(s),\bm{\gamma}(s)\rangle\bm{t}(s)+\mu \bm{n}(s),\bm{\gamma}(s)\rangle
=r\langle \bm{t}(s),\bm{\gamma}(s)\rangle^2+\mu\langle \bm{n}(s),\bm{\gamma}(s)\rangle.
\]
Suppose $\langle \bm{n}(s),\bm{\gamma}(s)\rangle\not =0.$ Then $F(s,\bm{x})=\partial H/\partial s(s,\bm{x})=0$ if and only if
\begin{eqnarray*}
\bm{x}&=& r\bm{\gamma}(s)+r\langle \bm{t}(s),\bm{\gamma}(s)\rangle \bm{t}(s)-\frac{r\langle \bm{t}(s),\bm{\gamma}(s)\rangle ^2 }{\langle \bm{n}(s),\bm{\gamma}(s)\rangle}\bm{n}(s) \\
&=&r\left(\bm{\gamma}(s)+\frac{\langle\bm{t}(s),\bm{\gamma}(s)\rangle\langle\bm{n}(s),\bm{\gamma}(s)\rangle\bm{t}(s)-\langle \bm{t}(s),\bm{\gamma}(s)\rangle ^2\bm{n}(s)}{\langle\bm{n}(s),\bm{\gamma}(s)
\rangle}\right).
\end{eqnarray*}
Since $\bm{\gamma}(s)=\langle \bm{\gamma}(s),\bm{t}(s)\rangle \bm{t}(s)+\langle \bm{\gamma}(s),\bm{n}(s)\rangle \bm{n}(s),$ the last expression equals to
\begin{eqnarray*}
&{}& r\left(\bm{\gamma}(s)+\frac{\langle\bm{n}(s),\bm{\gamma}(s)\rangle\bm{\gamma}(s)-(\langle\bm{n}(s),\bm{\gamma}(s)\rangle^2+\langle \bm{t}(s),\bm{\gamma}(s)\rangle ^2)\bm{n}(s)}{\langle\bm{n}(s),\bm{\gamma}(s)
\rangle}\right) \\
&{}& = r\left(2\bm{\gamma}(s)-\frac{(\langle\bm{n}(s),\bm{\gamma}(s)\rangle^2+\langle \bm{t}(s),\bm{\gamma}(s)\rangle ^2)}{\langle\bm{n}(s),\bm{\gamma}(s)
\rangle}\bm{n}(s)\right) \\
&{}&=r\left(2\bm{\gamma}(s)-\frac{\|\bm{\gamma}(s)\|^2}{\langle\bm{n}(s),\bm{\gamma}(s)
\rangle}\bm{n}(s)\right). 
\end{eqnarray*}
Therefore, the $r$-parallel primitivoid is parametrized by
\[
r\left(2\bm{\gamma}(s)-\frac{\|\bm{\gamma}(s)\|^2}{\langle\bm{n}(s),\bm{\gamma}(s)
\rangle}\bm{n}(s)\right),
\]
which is denoted by $r$-${\rm Pr}_{\bm{\gamma}}(s).$
Since 
\[
{\rm Pr}_{\bm{\gamma}}(s)=2\bm{\gamma}(s)-\frac{\|\bm{\gamma}(s)\|^2}{\langle\bm{n}(s),\bm{\gamma}(s)
\rangle}\bm{n}(s),
\]
we have the following theorem.
\begin{Th}
Suppose $\langle \bm{n}(s),\bm{\gamma}(s)\rangle\not =0.$ Then
we have $r$-${\rm Pr}_{\bm{\gamma}}(s)=r{\rm Pr}_{\bm{\gamma}}(s)={\rm Pr}_{r\bm{\gamma}}(s).$
\end{Th}
\demo
For $r\bm{\gamma}(s),$ the unit tangent vector is $r\bm{t}(s)/\|r\bm{t}(s)\|=\bm{t}(s),$
so that $\bm{n}(s)$ is the unit normal vector of $r\bm{\gamma}(s).$
Thus, we have
\[
{\rm Pr}_{r\bm{\gamma}}(s)=2r\bm{\gamma}(s)-\frac{\|r\bm{\gamma}(s)\|^2}{\langle\bm{n}(s),r\bm{\gamma}(s)
\rangle}\bm{n}(s)
=r\left(2\bm{\gamma}(s)-\frac{\|\bm{\gamma}(s)\|^2}{\langle\bm{n}(s),\bm{\gamma}(s)
\rangle}\bm{n}(s)\right).
\]
This completes the proof.
\enD

\section{Slant primitivoids}
Following the definition of the evolutoids of curves in the plane \cite{G-W}, we introduced 
the pedaloids of curves in \cite{I-T}.
For $\psi \in \R,$ define
\[
{\rm Pe}[\psi]_{\bm{\gamma}}(s)=\langle \bm{\gamma}(s),\cos\psi \bm{t}(s)+\sin\psi \bm{n}(s)\rangle (\cos\psi \bm{t}(s)+\sin\psi \bm{n}(s)).
\]
We call it a {\it $\psi$-pedaloid} of $\bm{\gamma}.$
If $\psi=\pi/2+n\pi, 3\pi/2+n\pi$ for $n\in \mathbb{Z},$ then ${\rm Pe}[\psi]_{\bm{\gamma}}(s)={\rm Pe}_{\bm{\gamma}}(s).$
If $\psi =n\pi$ for $n\in \mathbb{Z},$ then ${\rm Pe}[\psi]_{\bm{\gamma}}(s)=\langle \bm{\gamma}(s),\bm{t}(s)\rangle \bm{t}(s)$,
which is known as a {\it contrapedal} (i.e. a {\it C-pedal} for short) (cf. \cite{ZW}).
The C-pedal is denoted by ${\rm CPe}_{\bm{\gamma}}(s).$
Analogous to the notion of pedaloids, we define a $\phi$-slant primitivoid as the envelope of the family of lines
with the constant angle $\phi$ to the position vector of the curve.
Therefore, if $\phi=\pi/2+n\pi$ for $n\in \mathbb{Z},$ it is the primitive of the curve. We now give the precise definition
as follows: For a unit speed curve $\bm{\gamma}:I\lon \R^2\setminus \{\bm{0}\},$
we have $\bm{\gamma}(s)=\langle \bm{\gamma}(s),\bm{t}(s)\rangle \bm{t}(s)+\langle \bm{\gamma}(s),\bm{n}(s)\rangle \bm{n}(s),$ so that $\pi/2$-couter clockwise rotated vector is $J\bm{\gamma}(s)=-\langle \bm{\gamma}(s),\bm{n}(s)\rangle \bm{t}(s)+\langle \bm{\gamma}(s),\bm{t}(s)\rangle \bm{n}(s).$
For $\phi\in \R,$ we define 
\begin{eqnarray*}
\mathbb{N}[\phi](s)&=&\cos\phi \bm{\gamma}(s)+\sin\phi J\bm{\gamma}(s) \\
&=&\langle \bm{\gamma}(s),\cos\phi\bm{t}(s)-\sin\phi \bm{n}(s)\rangle \bm{t}(s)
+\langle \bm{\gamma}(s),\sin\phi\bm{t}(s)+\cos\phi \bm{n}(s)\rangle \bm{n}(s).
\end{eqnarray*}
We now consider a function $F:I\times \R^2\setminus\{\bm{0}\}\lon \R$ defined by
$F(s,\bm{x})=\langle \bm{x}-\bm{\gamma}(s),\mathbb{N}[\phi](s)\rangle.$
For any $s\in I,$ $f_s(\bm{x})=F(s,\bm{x})=0$ is a equation of the
line through $\bm{\gamma}(s)$ orthogonal to $\mathbb{N}[\phi](s)$, which means that
the angle between the line and the position vector $\bm{\gamma}(s)$ is $\phi+\pi/2.$
The envelope of the family of lines $\{f_s^{-1}(0)\}_{s\in I}$ is called 
a {\it $\phi$-slant primitivoid} of $\bm{\gamma}.$
By definition, the $n\pi$-primitivoid is the primitive of $\bm{\gamma}.$
We consider a parametrization of the $\phi$-slant primitivoid of $\bm{\gamma}.$
The $\phi$-slant primitivoid of $\bm{\gamma}$ is denoted by ${\rm Pr}[\phi]_{\bm{\gamma}}(s).$
\begin{Th} Let $\bm{\gamma}:I\lon \R^2\setminus \{\bm{0}\}$ be a unit speed
curve such that $\langle \bm{\gamma}(s),\bm{n}(s)\rangle\not= 0.$ Then we have
\[
{\rm Pr}[\phi]_{\bm{\gamma}}(s)=\cos\phi \left(\cos\phi{\rm Pr}_{\bm{\gamma}}(s)+\sin
\phi{\rm Pr}_{J\bm{\gamma}}(s)\right).
\] \end{Th}
\demo
Since $\mathbb{N}[\phi](s)=\cos\phi \bm{\gamma}(s)+\sin\phi J\bm{\gamma}(s),$
we have
$
\mathbb{N}[\phi]'(s)=\cos\phi \bm{t}(s)+\sin\phi J\bm{t}(s)=\cos\phi \bm{t}(s)+\sin\phi \bm{n}(s).
$
Therefore, $F(s,\bm{x})=0$ if and only if
there exits $\lambda \in \R$ such that 
$\bm{x}-\bm{\gamma}(s)=\lambda (\cos\phi J\bm{\gamma}(s)-\sin\phi\bm{\gamma}(s)).$
We remark that
$\langle J\bm{\gamma},\bm{n}\rangle =\langle \bm{\gamma},^t\!\!J\bm{n}\rangle
=-\langle \bm{\gamma},J\bm{n}\rangle=\langle \bm{\gamma},\bm{t}\rangle$
and
$\langle J\bm{\gamma},\bm{t}\rangle=\langle \bm{\gamma},^t\!\!J\bm{n}\rangle=
-\langle \bm{\gamma},J\bm{t}\rangle=-\langle\bm{\gamma},\bm{n}\rangle.$
With the condition $F(s,\bm{x})=0,$ we have  
\[
\frac{\partial F}{\partial s}(s,\bm{x})=-\langle \bm{t},\mathbb{N}[\phi]\rangle+\langle \bm{x}-\bm{\gamma},\mathbb{N}[\phi]'\rangle
=-(\cos\phi\langle\bm{t},\bm{\gamma}\rangle-\sin\phi \langle\bm{n},\bm{\gamma}\rangle) -\lambda \langle \bm{\gamma},\bm{n}\rangle.
\]
Therefore, $F(s,\bm{x})=\partial F(s,\bm{x})/\partial s=0$ if and only if
\[
\bm{x}-\bm{\gamma}(s) =-\frac{\langle \bm{\gamma}(s),\cos\phi\bm{t}(s)-\sin\phi \bm{n}(s)\rangle}{\langle\bm{\gamma}(s),\bm{n}(s)\rangle}(\cos\phi J\bm{\gamma}(s)-\sin\phi\bm{\gamma}(s)).
\]
Since $\bm{\gamma}=\langle \bm{\gamma},\bm{t}\rangle \bm{t}+\langle\bm{\gamma},\bm{n}\rangle\bm{n},$
we have
$J\bm{\gamma}=-\langle\bm{\gamma},\bm{n}\rangle\bm{t}+\langle \bm{\gamma},\bm{t}\rangle \bm{n}.$
Moreover, we have $\|\bm{\gamma}(s)\|^2=\langle \bm{\gamma}(s),\bm{t}(s)\rangle^2+\langle \bm{\gamma}(s),\bm{n}(s)\rangle^2,$ so that
\begin{eqnarray*}
\bm{x}&{}&=-\frac{\langle \bm{\gamma}(s),\cos\phi\bm{t}(s)-\sin\phi \bm{n}(s)\rangle}{\langle\bm{\gamma}(s),\bm{n}(s)\rangle}(\cos\phi J\bm{\gamma}(s)-\sin\phi\bm{\gamma}(s))+\bm{\gamma}(s)  \\
&{}&= \frac{\cos\phi}{\langle \bm{\gamma}(s),\bm{n}(s)\rangle}((-\cos\phi \langle \bm{\gamma}(s),\bm{t}(s)\rangle +\sin\phi \langle \bm{\gamma}(s),\bm{n}(s)\rangle)J\bm{\gamma}(s)\\
&{}&\qquad\qquad +(\sin\phi \langle \bm{\gamma}(s),\bm{t}(s)\rangle+\cos\phi \langle \bm{\gamma}(s),\bm{n}(s)\rangle)\bm{\gamma}(s)) \\
&{}&= \frac{\cos\phi}{\langle \bm{\gamma}(s),\bm{n}(s)\rangle}((2\cos\phi \langle \bm{\gamma}(s),\bm{t}(s)\rangle\langle \bm{\gamma}(s),\bm{n}(s)\rangle+\sin\phi (\langle \bm{\gamma}(s),\bm{t}(s)\rangle^2-\langle \bm{\gamma}(s),\bm{n}(s)\rangle^2)\bm{t}(s) \\
&{}& \qquad\qquad + (\cos\phi (\langle \bm{\gamma}(s),\bm{n}(s)\rangle^2-\langle \bm{\gamma}(s),\bm{t}(s)\rangle^2)+
2\sin\phi \langle \bm{\gamma}(s),\bm{t}(s)\rangle\langle \bm{\gamma}(s),\bm{n}(s)\rangle)\bm{n}(s)) \\
&{}& =\frac{\cos\phi}{\langle \bm{\gamma}(s),\bm{n}(s)\rangle}(2\langle \bm{\gamma}(s),\bm{t}(s)\rangle
\langle \bm{\gamma}(s),\bm{n}(s)\rangle(\cos\phi \bm{t}(s)+\sin\phi \bm{n}(s))\\
&{}& \qquad \qquad +\sin\phi (2\langle \bm{\gamma}(s),\bm{t}(s)\rangle^2-\|\bm{\gamma}(s)\|^2)\bm{t}(s)+\cos\phi (2\langle \bm{\gamma}(s),\bm{n}(s)\rangle^2-\|\bm{\gamma}(s)\|^2)\bm{n}(s)) \\
&{}& = \frac{\cos\phi}{\langle \bm{\gamma}(s),\bm{n}(s)\rangle}(\cos\phi (2\langle \bm{\gamma}(s),\bm{n}(s)\rangle\bm{\gamma}(s) -\|\bm{\gamma}(s)\|^2\bm{n}(s)) \\
&{}& \qquad \qquad \hskip 5cm +\sin\phi (2\langle \bm{\gamma}(s),\bm{t}(s)\rangle\bm{\gamma}(s) -\|\bm{\gamma}(s)\|^2\bm{t}(s)) \\
&{}&=\cos\phi \left(\cos\phi \left(2\bm{\gamma}(s)-\frac{\|\bm{\gamma}(s)\|^2}{\langle
\bm{\gamma}(s),\bm{n}\rangle}\bm{n}(s)\right)+\sin\phi \left(2\frac{\langle\bm{\gamma}(s),\bm{t}(s)\rangle}{
\langle\bm{\gamma}(s),\bm{n}(s)\rangle}\bm{\gamma}(s)-\frac{\|\bm{\gamma}(s)\|^2}{\langle
\bm{\gamma}(s),\bm{n}\rangle}\bm{t}(s)\right)\right).
\end{eqnarray*}
Here, ${\rm Pr}_{\bm{\gamma}}(s)=2\bm{\gamma}(s)-\frac{\|\bm{\gamma}(s)\|^2}{\langle
\bm{\gamma}(s),\bm{n}\rangle}\bm{n}(s),$ $J\bm{\gamma}(s)$ is a unit speed curve such that the tangent vector is $J\bm{t}(s)=\bm{n}(s)$ and the normal vector is $J\bm{n}(s)=JJ\bm{t}(s)=-\bm{t}(s).$
Therefore, we have
\[
{\rm Pr}_{J\bm{\gamma}}(s)=2J\bm{\gamma}(s)-\frac{\|J\bm{\gamma}(s)\|^2}{\langle J\bm{\gamma}(s),J\bm{n}(s)\rangle}J\bm{n}(s)=2J\bm{\gamma}(s)+\frac{\|\bm{\gamma}(s)\|^2}{\langle \bm{\gamma}(s),\bm{n}(s)\rangle}\bm{t}(s).
\]
Since $J\bm{\gamma}=-\langle\bm{\gamma},\bm{n}\rangle\bm{t}+\langle \bm{\gamma},\bm{t}\rangle \bm{n}$ and 
$\bm{\gamma}=\langle \bm{\gamma},\bm{t}\rangle \bm{t}+\langle\bm{\gamma},\bm{n}\rangle\bm{n},$
the last formula is equal to
\begin{eqnarray*}
&{}&2(-\langle\bm{\gamma}(s),\bm{n}(s)\rangle\bm{t}(s)+\langle \bm{\gamma}(s),\bm{t}(s)\rangle \bm{n}(s))+\frac{\|\bm{\gamma}(s)\|^2}{\langle \bm{\gamma}(s),\bm{n}(s)\rangle}\bm{t}(s) \\
&{}&=\frac{2}{\langle \bm{\gamma}(s),\bm{n}(s)\rangle }(-\langle\bm{\gamma}(s),\bm{n}(s)\rangle^2\bm{t}(s)+\langle\bm{\gamma}(s),\bm{n}(s)\rangle\langle\bm{\gamma}(s),\bm{t}(s)\rangle \bm{n}(s))+\frac{\|\bm{\gamma}(s)\|^2}{\langle \bm{\gamma}(s),\bm{n}(s)\rangle}\bm{t}(s) \\
&{}& =\frac{2}{\langle \bm{\gamma}(s),\bm{n}(s)\rangle }(\langle\bm{\gamma}(s),\bm{t}(s)\rangle(\bm{\gamma}(s)-\langle \bm{\gamma}(s),\bm{t}(s)\rangle\bm{t}(s)-\langle\bm{\gamma}(s),\bm{n}(s)\rangle^2\bm{t}(s)))+\frac{\|\bm{\gamma}(s)\|^2}{\langle \bm{\gamma}(s),\bm{n}(s)\rangle}\bm{t}(s) \\
&{}&=\frac{2}{\langle \bm{\gamma}(s),\bm{n}(s)\rangle }(\langle\bm{\gamma}(s),\bm{t}(s)\rangle\bm{\gamma}(s)-\|\bm{\gamma}(s)\|^2\bm{t}(s))+\frac{\|\bm{\gamma}(s)\|^2}{\langle \bm{\gamma}(s),\bm{n}(s)\rangle}\bm{t}(s)
\\
&{}&=2\frac{\langle\bm{\gamma}(s),\bm{t}(s)\rangle}{
\langle\bm{\gamma}(s),\bm{n}(s)\rangle}\bm{\gamma}(s)-\frac{\|\bm{\gamma}(s)\|^2}{\langle
\bm{\gamma}(s),\bm{n}\rangle}\bm{t}(s).
\end{eqnarray*}
Thus we have
\[
{\rm Pr}_{J\bm{\gamma}}(s)=2\frac{\langle\bm{\gamma}(s),\bm{t}(s)\rangle}{
\langle\bm{\gamma}(s),\bm{n}(s)\rangle}\bm{\gamma}(s)-\frac{\|\bm{\gamma}(s)\|^2}{\langle
\bm{\gamma}(s),\bm{n}\rangle}\bm{t}(s).
\]
This completes the proof. \enD
\begin{Exa}\rm
We consider the ellipse defined by $x_1^2+3x_2^2=1.$
By Theorem 4.1, the $\phi$-slant primitivoid is parametrized by
\begin{eqnarray*}
{\rm Pr}[\phi]_{\bm{\gamma}}(t)&=&\left(\frac{\cos\phi}{6}\left(-2\cos t(-4+\cos 2t)\cos\phi+\sqrt{3}(-\sin t+\sin 3t)\sin\phi\right),\right. \\
&{}& \qquad \left.\left(\frac{\cos\phi}{3}(4\cos t\sin\phi-\cos 2t (\sqrt{3}\cos\phi\sin t+\cos t\sin\phi))\right)\right).
\end{eqnarray*}

\begin{center}
\begin{tabular}{cc}   
\includegraphics[height=50mm,angle=90]{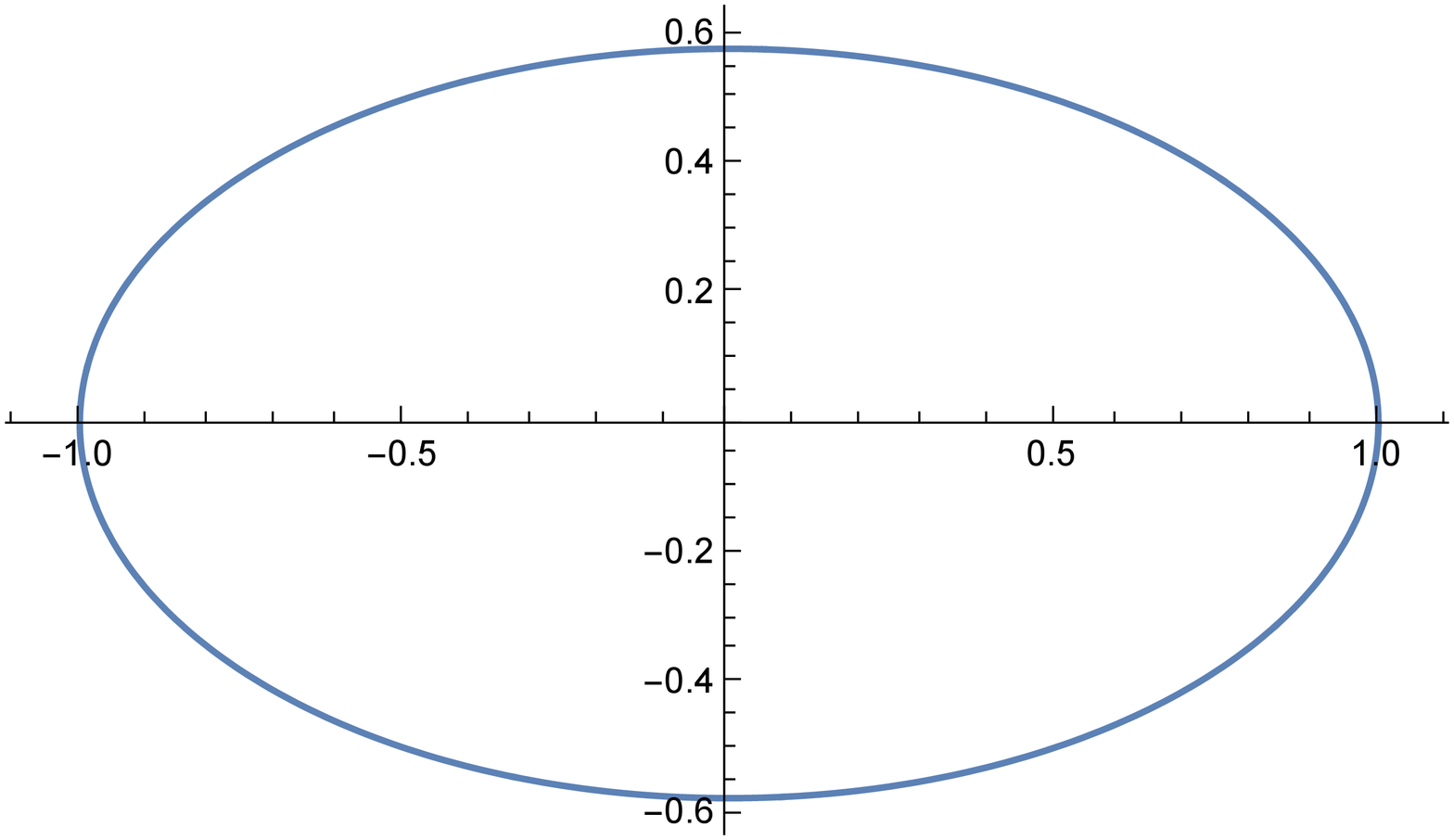}
&
\includegraphics[width=35mm,angle=90]{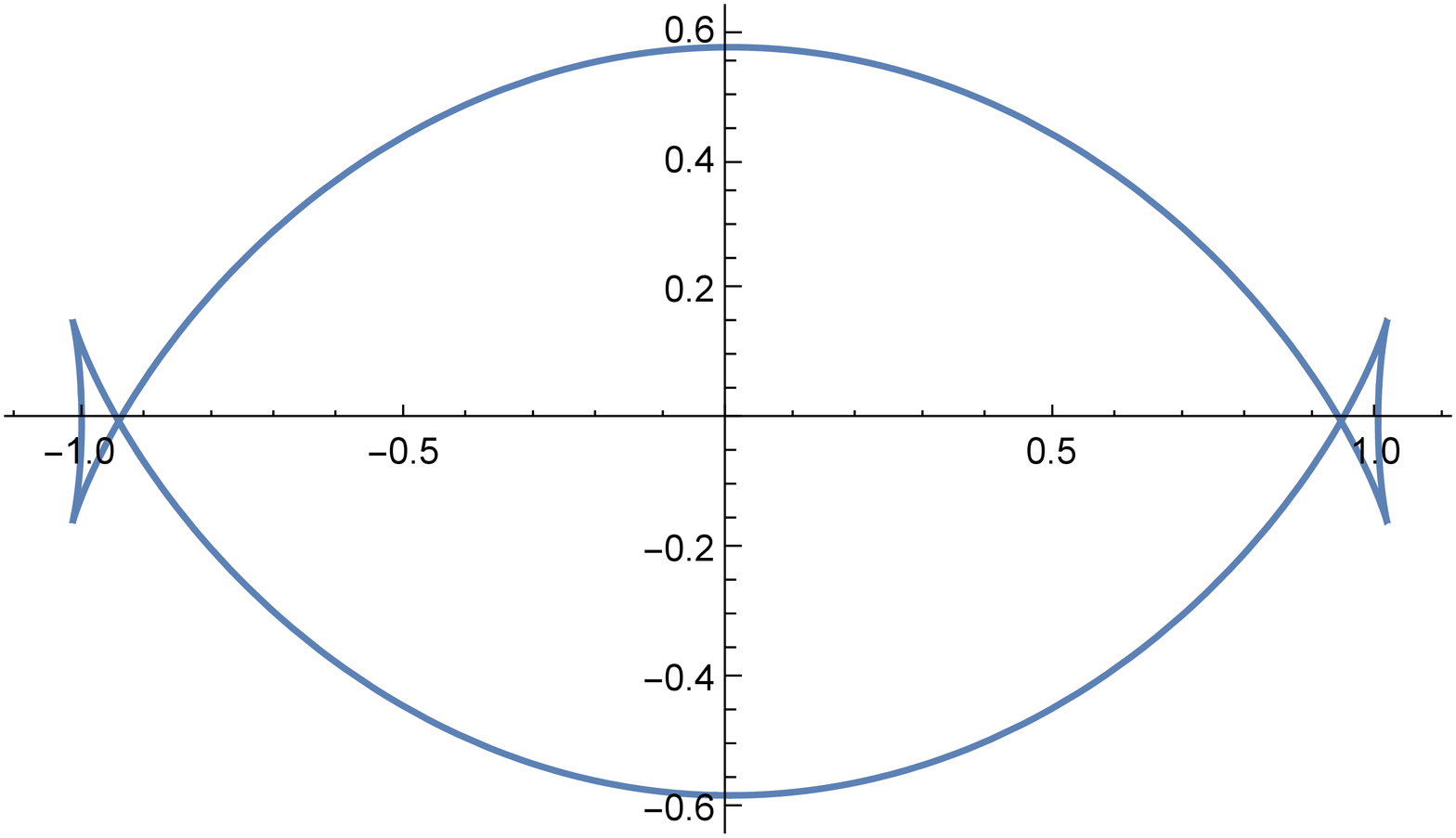}
 \\
\text{Fig. 1 : The ellipse}
&
\text{Fig. 2 : The 0-primitivoid=the primitive}
\end{tabular}
\end{center}

\begin{center}
\begin{tabular}{cc}   
\includegraphics[height=40mm,angle=90]{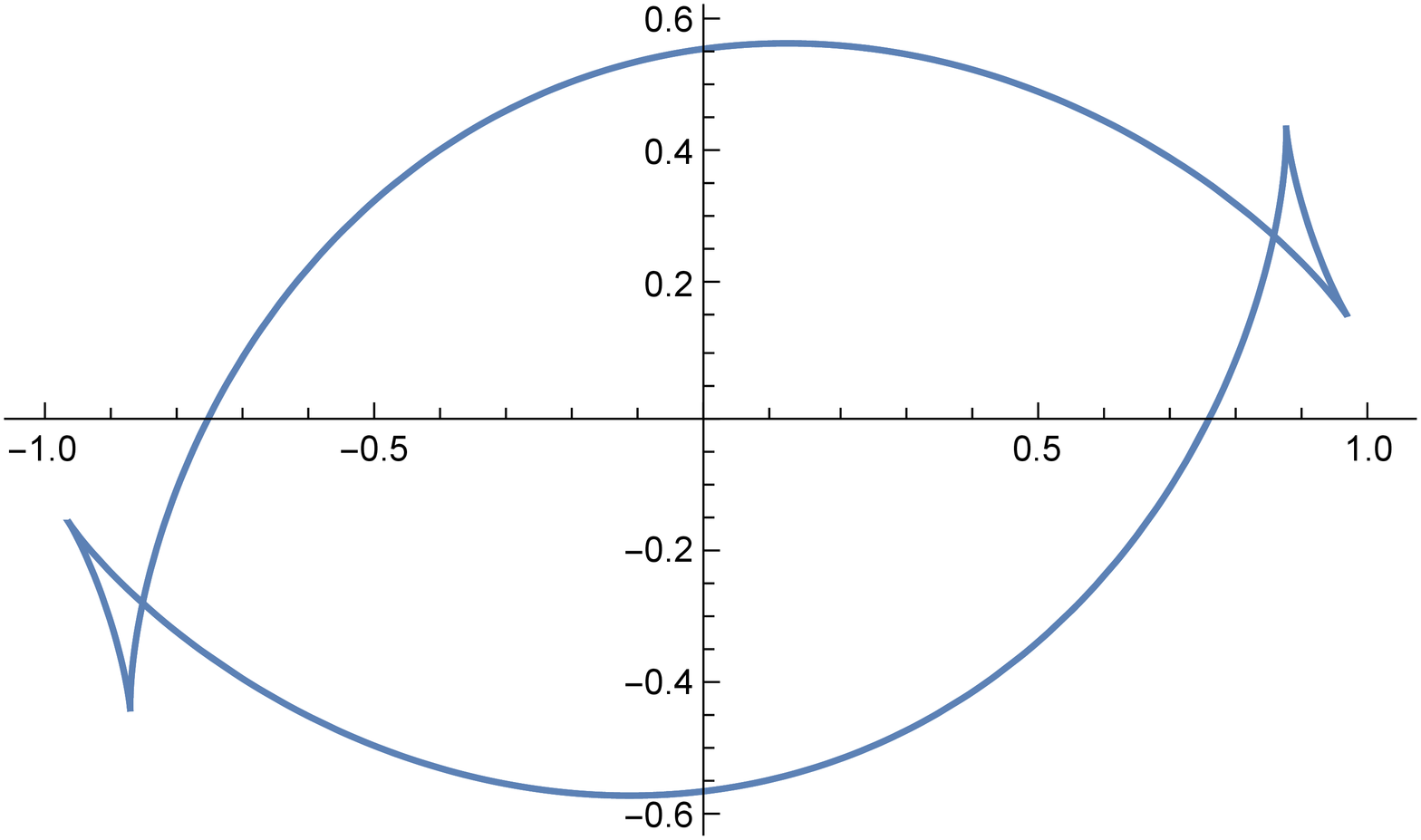}
&
\includegraphics[width=28mm,angle=90]{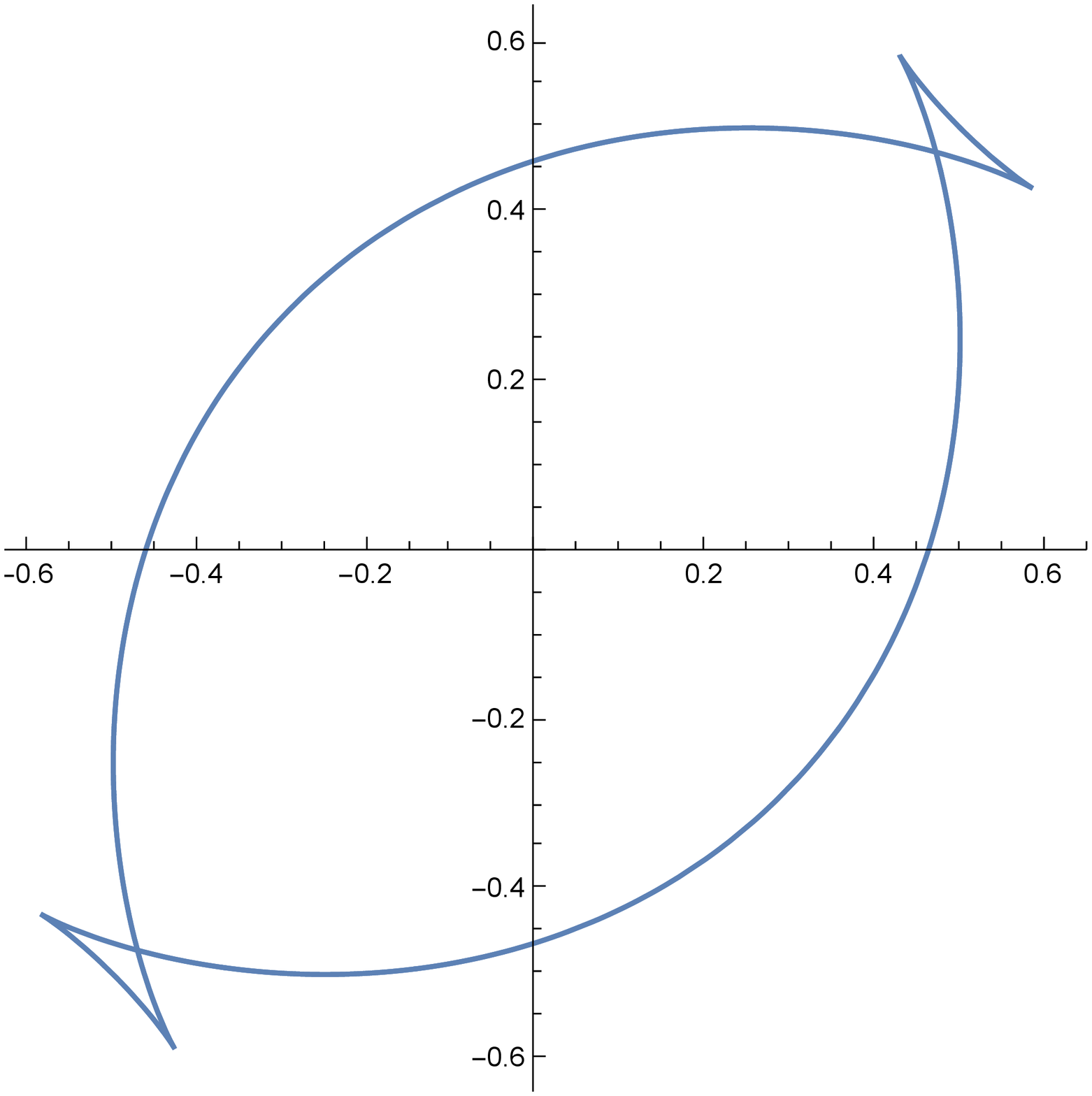}
 \\
\text{Fig. 3 : The $\pi/10$-primitivoid}
&
\text{Fig. 4 : The $\pi/4$-primitivoid}
\end{tabular}
\end{center}

\begin{center}
\begin{tabular}{cc}   
\includegraphics[height=28mm,angle=90]{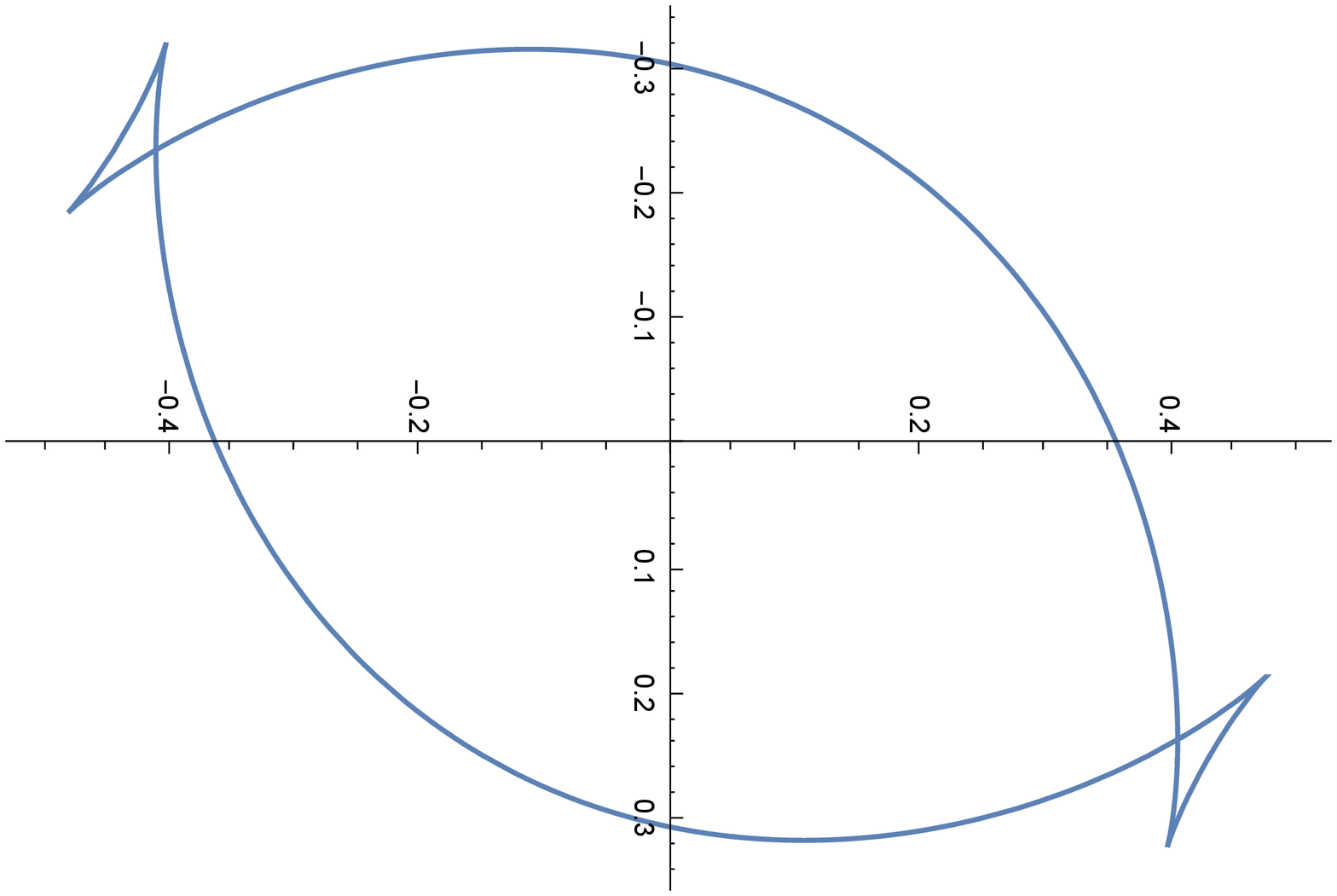}
&
\includegraphics[width=30mm,angle=90]{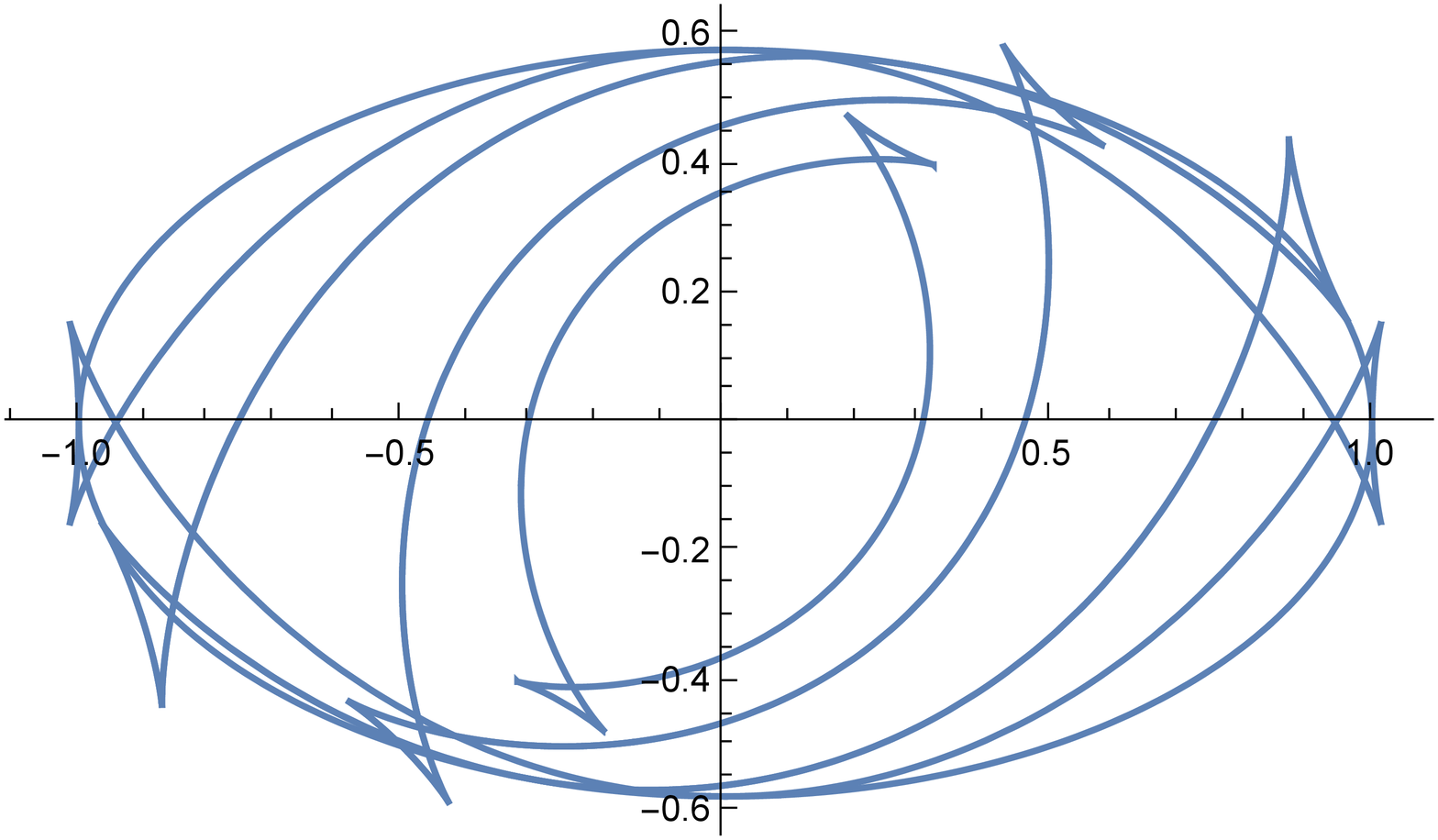}
 \\
\text{Fig. 5 : The $\pi/3$-primitivoid}
&
\text{Fig. 6 : The ellipse and primitivoids}
\end{tabular}
\end{center}
We observe in Fig. 6 that these primitivoids might be similar in shape.
\end{Exa}
\par
For $\phi\in [0,2\pi),$ we denote the rotation matrix by
\[
R(\phi)=\begin{pmatrix} \cos\phi & -\sin\phi \\ \sin\phi & \cos\phi \end{pmatrix}.
\]
Let $\bm{\gamma}:I\lon \R^2\setminus \{\bm{0}\}$ be a unit speed curve.
We denote the {\it $\phi$-rotated curve} by
$\bm{\gamma}_{\phi} (s)=R(\phi)\bm{\gamma}(s).$
Since $\phi$-rotated curve is also a unit speed curve, we have $\bm{t}_{\phi}(s)=\bm{\gamma}_{\phi}' (s)=R(\phi)\bm{t}(s)$ and
$\bm{n}_{\phi}(s)=J\bm{t}_{\phi}(s)=JR(\phi)\bm{t}(s)=R(\phi)J\bm{t}(s)=R(\phi)\bm{n}(s).$ 
It follows that 
\[
{\rm Pr}_{\bm{\gamma}_{\phi}}(s)=R(\phi){\rm Pr}_{\bm{\gamma}}(s).
\]
Moreover, we have the following lemma.
\begin{Lem} With the same notations as above, we have
\[
J{\rm Pr}_{\bm{\gamma}}(s)={\rm Pr}_{J\bm{\gamma}}(s).
\]
\end{Lem}
\demo
We have
\begin{eqnarray*}
J{\rm Pr}_{\bm{\gamma}}(s)&=&
J\left(2\bm{\gamma}(s)-\frac{\|\bm{\gamma}(s)\|^2}{\langle
\bm{\gamma}(s),\bm{n}(s)\rangle}\bm{n}(s)\right)\\
&=&2J\bm{\gamma}(s)-\frac{\|\bm{\gamma}(s)\|^2}{\langle
\bm{\gamma}(s),\bm{n}(s)\rangle}J\bm{n}(s)=2J\bm{\gamma}(s)+\frac{\|\bm{\gamma}(s)\|^2}{\langle
\bm{\gamma}(s),\bm{n}(s)\rangle}\bm{t}(s).
\end{eqnarray*}
In the proof of Theorem 4.1, we have shown that
\[
{\rm Pr}_{J\bm{\gamma}}(s)=2J\bm{\gamma}(s)+\frac{\|\bm{\gamma}(s)\|^2}{\langle
\bm{\gamma}(s),\bm{n}(s)\rangle}\bm{t}(s).
\]
This completes the proof. \enD
Then we have the following theorem as a corollary of Theorem 4.1.
\begin{Th} With the same notations as above, we have
\[
{\rm Pr}[\phi]_{\bm{\gamma}}(s)=\cos\phi R(\phi){\rm Pr}_{\bm{\gamma}}(s)=\cos\phi {\rm Pr}_{\bm{\gamma}_{\phi}}(s)={\rm Pr}_{\cos\phi\bm{\gamma}_{\phi}}(s).
\]
\end {Th}
\demo
It follows from Theorem 4.1 and Lemma 4.3 that
\[
{\rm Pr}[\phi]_{\bm{\gamma}}(s)=\cos\phi \left(
\cos\phi {\rm Pr}_{\bm{\gamma}}(s)+\sin\phi {\rm Pr}_{\bm{J\gamma}}(s)
\right)=\cos\phi \left(
\cos\phi {\rm Pr}_{\bm{\gamma}}(s)+\sin\phi J{\rm Pr}_{\bm{\gamma}}(s)
\right).
\]
Moreover, for $\bm{a}\in \R^2\setminus \{\bm{0}\}$, if we set $\bm{b}=J\bm{a},$
then we have
\[
\cos\phi \bm{a}+\sin\phi \bm{b}= \cos\phi \bm{a}+\sin\phi J\bm{a}
=(\cos\phi I+\sin\phi J)\bm{a}=R(\phi)\bm{a}.
\]
Therefore, we have
\[
{\rm Pr}[\phi]_{\bm{\gamma}}(s)=\cos\phi R(\phi){\rm Pr}_{\bm{\gamma}}(s).
\]
This completes the proof. \enD
\par
We now consider the relation of primitivoids with pedals.
For $\lambda \in \R\setminus \{0\}, $ we consider $\lambda \bm{\gamma}(s).$
Since $(\lambda\bm{\gamma})'(s)=\lambda \bm{t}(s),$ we have
$\|\lambda\bm{t}(s)\|=|\lambda|$, so that
$\bm{t}_{\lambda\bm{\gamma}}(s)=\bm{t}(s)$ for $\lambda >0$
and $\bm{n}_{\lambda\bm{\gamma}}(s)=J\bm{t}(s)=\bm{n}(s).$
We also have $\bm{t}_{\lambda\bm{\gamma}}(s)=-\bm{t}(s)$ for $\lambda <0$
and $\bm{n}_{\lambda\bm{\gamma}}(s)=J(-\bm{t}(s))=-\bm{n}(s).$
Thus, we have
${\rm Pe}_{\lambda\bm{\gamma}}(s)=\lambda {\rm Pe}_{\bm{\gamma}}(s).$
Then we have the following proposition.
\begin{Pro}
With the same notations as above, suppose that ${\rm Pe}_{\bm{\gamma}}$ and ${\rm Pr}[\phi]_{\bm{\gamma}}$
are regular curves. Then we have
\[
{\rm Pr}[\phi]_{{\rm Pe}_{\bm{\gamma}}}(s)={\rm Pe}_{{\rm Pr}[\phi]_{\bm{\gamma}}}(s)
=\cos\phi R(\phi)\bm{\gamma}(s)=\cos\phi \bm{\gamma}_\phi (s).
\]
\end{Pro}
\demo
By the above arguments and Theorem 4.4, we have
\[
{\rm Pe}_{{\rm Pr}[\phi]_{\bm{\gamma}}}(s)={\rm Pe}_{\cos\phi {\rm Pr}_{\bm{\gamma}_\phi}}(s)
=\cos\phi {\rm Pe}_{{\rm Pr}_{\bm{\gamma}_\phi}}(s)=\cos\phi \bm{\gamma}_\phi(s)
\]
and
\[
{\rm Pr}[\phi]_{{\rm Pe}_{\bm{\gamma}}}(s)=\cos\phi R(\phi) {\rm Pr}_{{\rm Pe}_{\bm{\gamma}}}(s)
=\cos\phi R(\phi) \bm{\gamma}(s)=\cos\phi \bm{\gamma}_\phi (s).
\]
This completes the proof. \enD
\par
We now consider relations of primitivoids with anti-pedals and parallel primitivoids.
As a corollary of Theorem 4.4, we have the following proposition.
\begin{Pro} Let $\bm{\gamma}:I\lon \R^2\setminus \{\bm{0}\}$ be a unit speed curve
such that $\langle \bm{n}(s),\bm{\gamma}(s)\rangle \not=0.$
Then we have \[
{\rm Pr}[\phi]_{\bm{\gamma}}(s)=\cos\phi {\rm APe}_{\Psi\circ \bm{\gamma}_\phi}(s)
=R(\phi)(\cos\phi\mbox{\rm-Pr}_{\bm{\gamma}})(s).\]
\end{Pro}
\demo
By Proposition 2.1 and Theorem 4.4, we have
\[
{\rm Pr}[\phi]_{\bm{\gamma}}(s)=\cos\phi {\rm Pr}_{\bm{\gamma}_{\phi}}(s)
=\cos\phi {\rm APe}_{\Psi\circ\bm{\gamma}_\phi}(s).
\]
By Theorem 3.1, we have
\[
\cos\phi {\rm Pr}_{\bm{\gamma}_{\phi}}(s)=\cos\phi R(\phi){\rm Pr}_{\bm{\gamma}}(s)
=R(\phi)\cos\phi {\rm Pr}_{\bm{\gamma}}(s)=R(\phi)(\cos\phi\mbox{\rm-Pr}_{\bm{\gamma}})(s).
\]
\enD

\section{Pedals and primitivoids of frontals}
In the previous sections we investigated pedals and primitivoids of regular curves.
However, pedals and primitivoids generally have singularities even for regular curves.
In the last section we consider the pedal (respectively, the primitivoid) of the primitivoid (respectively, the pedal) of a curve.
Therefore, we need to generalize the notions of primitivoids (respectively, pedals) of certain singular curves.
One of the natural singular curves in the Euclidean plane for which we can develop the differential geometry is the notion of frontals \cite{F-T, F-T2}.
\par
We say that $(\bm{\gamma},\bm{\nu}):I\lon \R^{2}\times S^{1}$ is a {\it Legendrian curve}
if $(\bm{\gamma},\bm{\nu})^{*}\theta =0$, where $\theta$ is the canonical contact
$1$-form on the unit tangent bundle $T_{1}\R^{2}=\R^{2}\times S^{1}$ (cf. \cite{Arnold2}). 
The last condition is equivalent to $\langle \dot{\bm{\gamma}}(t),\bm{\nu}(t)\rangle=0$
for any $t\in I.$
We say that $\bm{\gamma}:I\lon \R^{2}$ is a {\it frontal} if
there exists $\bm{\nu}:I\lon S^{1}$ such that $(\bm{\gamma},\bm{\nu})$ is a Legendrian curve.
If $(\bm{\gamma},\bm{\nu})$ is an immersion, $\bm{\gamma}$ is said to be a {\it front}.
A differential geometry on frontals was constructed in \cite{F-T2}.
For a Legendrian curve $(\bm{\gamma},\bm{\nu}):I\lon \R^{2}\times S^{1},$
we define a unit vector field $\bm{\mu}(t)=J(\bm{\nu}(t))$ along $\bm{\gamma}.$ 
Then we have the following Frenet type formulae \cite{F-T2}:
\[
\left\{
\begin{array}{ll}
\dot{\bm{\nu}}(t)=\ell(t)\bm{\mu}(t), \\
\dot{\bm{\mu}}(t)=-\ell (t)\bm{\nu}(t),
\end{array}
\right.
\]
where $\ell (t)=\langle \dot{\bm{\nu}}(t),\bm{\mu}(t)\rangle.$
Moreover, there exists $\beta (t)$ such that $\dot{\bm{\gamma}}(t)=\beta(t)\bm{\mu}(t)$
for any $t\in I.$
The pair $(\ell,\beta)$ is called a {\it curvature of the Legendrian curve} $(\bm{\gamma},\bm{\nu}).$
By definition, $t_{0}\in I$ is a singular point of $\bm{\gamma}$ if and only if 
$\beta (t_{0})=0.$ Moreover, for a regular curve $\bm{\gamma}$, $\bm{\mu}(t)=\bm{t}(t)$ and 
$\ell (t)=\|\dot{\bm{\gamma}}(t)\|\kappa (t).$
The Legendrian curve $(\bm{\gamma},\bm{\nu})$ is immersive (i.e. $\bm{\gamma}$ is a front) if and only if 
$(\ell(t),\beta(t))\not= (0,0)$ for any $t\in I .$ 
So the inflection point $t_{0}\in I$ of the frontal $\bm{\gamma}$ is a point $\ell(t_{0})=0.$
For more detailed properties of Legendrian curves, see \cite{F-T, F-T2}.
\par
In \cite{I-T,LP}, the {\it pedal} of a frontal $\bm{\gamma}$ is defined by
\[
\mathcal{P}e_{\bm{\gamma}}(t)=\langle \bm{\gamma}(t), \bm{\nu}(t)\rangle \bm{\nu}(t).
\]
We have shown that if there exist $\delta (t)$ and 
$\bm{\sigma}:I\lon S^{1}$such that $\bm{\gamma} (t)=\delta (t)\bm{\sigma} (t)$
for any $t\in I,$  the pedal $\mathcal{P}e_{\bm{\gamma}}(t)$ of
$\bm{\gamma}$ is a frontal.
Moreover, the {\it anti-pedal} of $\bm{\gamma}$ is defined to be
\[
\mathcal{AP}e_{\bm{\gamma}}(s)=\frac{1}{\langle \bm{\gamma}(t), \bm{\nu}(t)\rangle} \bm{\nu}(t).
\]
As we remarked for a regular curve $\bm{\gamma}$ in \S 1,
the anti-pedal of a frontal $\bm{\gamma}$ is the envelope of the family of lines
$\{\bm{x}\ |\ \langle \bm{x},\bm{\gamma}(t)\rangle -1=0\}_{t\in I}.$

On the other hand, for a frontal $\bm{\gamma}$
with $\langle
\bm{\gamma}(t),\bm{\nu}(t)\rangle\not=0,$ we define the {\it primitive} of $\bm{\gamma}$
by
\[
\mathcal{P}r_{\bm{\gamma}}(t)=2\bm{\gamma}(t)-\frac{\|\bm{\gamma}(t)\|^2}{\langle
\bm{\gamma}(t),\bm{\nu}(t)\rangle}\bm{\nu}(t).
\]
We remark that the primitive of $\bm{\gamma}$ is the envelope of
the family of lines \[
\{\bm{x}\ |\ \langle \bm{x}-\bm{\gamma}(t),\bm{\gamma}(t)\rangle =0\}_{t\in I}.
\]
Since
$
\langle \bm{x}-\bm{\gamma}(t),\bm{\gamma}(t)\rangle =\langle \bm{x},\bm{\gamma}(t)\rangle-\|\bm{\gamma}(t)\|^2,
$
$\langle \bm{x}-\bm{\gamma}(t),\bm{\gamma}(t)\rangle =0$ if and only if
$\langle \bm{x},\Psi\circ \bm{\gamma}(t)\rangle-1=0,$
where $\Psi:\R^2\setminus\{\bm{0}\}\lon \R^2\setminus\{\bm{0}\}$ is the inversion
defined by $\Psi (\bm{x})=\frac{\bm{x}}{\|\bm{x}\|^2}.$
Therefore we have the following lemma.
\begin{Lem}
Let $(\bm{\gamma},\bm{\nu}):I\lon (\R^{2}\setminus\{\bm{0}\})\times S^{1}$ be a Legendrian curve
such that $\langle
\bm{\gamma}(t),\bm{\nu}(t)\rangle\not=0.$
Then we have
\[
\mathcal{P}r_{\bm{\gamma}}(t)=\mathcal{AP}e_{\Psi\circ\bm{\gamma}}(t)=\Psi\circ\mathcal{P}e_{\Psi\circ\bm{\gamma}}(t).
\]
\end{Lem}
Then we have the following lemma.
\begin{Lem} Let $(\bm{\gamma},\bm{\nu}):I\lon (\R^{2}\setminus\{\bm{0}\})\times S^{1}$ be a Legendrian curve
such that $\langle
\bm{\gamma}(t),\bm{\nu}(t)\rangle\not=0.$ Then the primitive $\mathcal{P}r_{\bm{\gamma}}$ of 
$\bm{\gamma}$ is a frontal.
\end{Lem}
\demo
Since $\bm{\gamma}(t)\not= \bm{0},$ $\Psi\circ \bm{\gamma}(t)$ is well defined and it is not 
equal to the origin. Thus, $\mathcal{P}e_{\Psi\circ\bm{\gamma}}$ is a frontal.
By Lemma 5.1, we have $\mathcal{P}r_{\bm{\gamma}}(t)=\Psi\circ\mathcal{P}e_{\Psi\circ\bm{\gamma}}(t).$
Here, $\Psi$ is a diffeomorphism, so that $\Psi\circ\mathcal{P}e_{\Psi\circ\bm{\gamma}}$ is a
frontal.
\enD
\par
Once we obtained parametrization of primitivoids of regular curves, we can generalize these notions for
frontals.
Following Theorems 3.1 and 4.4,
we respectively define the {\it $r$-parallel primitivoid} and the {\it $\phi$-slant primitivoid} of a frontal $\bm{\gamma}$ with $\langle
\bm{\gamma}(t),\bm{\nu}(t)\rangle\not=0$ by
\[
r\mbox{\rm -}\mathcal{P}r_{\bm {\gamma}}(t)=r\mathcal{P}r_{\bm{\gamma}}(t)\ \mbox{\rm and}\
\mathcal{P}r[\phi]_{\bm{\gamma}}(t)=\cos\phi R(\phi)\mathcal{P}r_{\bm{\gamma}}(t).
\]
Since $\langle R(\phi)\bm{\gamma}(t), R(\phi)\bm{\nu}(t)\rangle =\langle \bm{\gamma}(t),\bm{\nu}(t)\rangle,$
we have
\[
R(\phi)\mathcal{P}r_{\bm{\gamma}}(t)=R(\phi)\left(2\bm{\gamma}(t)-\frac{\|\bm{\gamma}(t)\|^2}{\langle
\bm{\gamma}(t),\bm{\nu}(t)\rangle}\bm{\nu}(t)\right)
=\mathcal{P}r_{R(\phi)\bm{\gamma}}(t).
\]
A {\it rotated frontal $\bm{\gamma}_\phi$} is defined to be $\bm{\gamma}_\phi(t)=R(\phi)\bm{\gamma}(t).$
Thus, we have $R(\phi)\mathcal{P}r_{\bm{\gamma}}(t)=\mathcal{P}r_{\bm{\gamma}_\phi}(t).$
Then we have the following theorem.
\begin{Th}
Let $(\bm{\gamma},\bm{\nu}):I\lon (\R^{2}\setminus\{\bm{0}\})\times S^{1}$ be a Legendrian curve
such that $\langle
\bm{\gamma}(t),\bm{\nu}(t)\rangle\not=0.$ Then both of the $r$-parallel primitivoid 
$r\mbox{\rm -}\mathcal{P}r_{\bm {\gamma}}(t)$ and the $\phi$-slant primitivoid $\mathcal{P}r[\phi]_{\bm{\gamma}}(t)$ of 
$\bm{\gamma}$ are frontals.
\end{Th}
\demo
By Lemma 5.2, the primitive $\mathcal{P}r_{\bm{\gamma}}$ is a frontal.
Since linear map $\cos\phi R(\phi):\R^2\lon \R^2$ is a diffeomorphism,
$\mathcal{P}r[\phi]_{\bm{\gamma}}=\cos\phi R(\phi)\mathcal{P}r_{\bm{\gamma}}$ is a frontal.  
Moreover, if there exists $t_0\in I$ such that
$\mathcal{P}r_{\bm{\gamma}}(t_0)=0,$ then we have
\[
2\bm{\gamma}(t_0)=\frac{\|\bm{\gamma}(t_0)\|^2}{\langle\bm{\gamma}(t_0),\bm{\nu}(t_0)\rangle}\bm{\nu}(t_0).
\]
Since $\|\bm{\nu}(t_0)\|=1,$ we have
$\|\bm{\gamma} (t_0)\|=2|\langle\bm{\gamma}(t_0),\bm{\nu}(t_0)\rangle|.$
It follows that
\[
2\bm{\gamma}(t_0)=\frac{4\langle\bm{\gamma}(t_0),\bm{\nu}(t_0)\rangle^2}{\langle\bm{\gamma}(t_0),\bm{\nu}(t_0)\rangle}\bm{\nu}(t_0),
\]
so that $\bm{\gamma}(t_0)=2\langle\bm{\gamma}(t_0),\bm{\nu}(t_0)\rangle\bm{\nu}(t_0).$
Thus, we have $\langle\bm{\gamma}(t_0),\bm{\nu}(t_0)\rangle=2\langle\bm{\gamma}(t_0),\bm{\nu}(t_0)\rangle.$
This contradicts to the assumption $\langle
\bm{\gamma}(t),\bm{\nu}(t)\rangle\not=0.$
Therefore, $\mathcal{P}r_{\bm{\gamma}}(t)\not=0.$
We now define $\psi_r:\R^2\setminus \{0\}\lon \R^2\setminus \{0\}$ by $\psi_r(\bm{x})=r\bm{x}.$
Then $\psi_r$ is a diffeomorphism and $\psi_r\circ \mathcal{P}r_{\bm{\gamma}}=r\mbox{\rm -}\mathcal{P}r_{\bm {\gamma}}.$
It follows that $r\mbox{\rm -}\mathcal{P}r_{\bm {\gamma}}$ is a frontal.
\enD
By the above theorem, we can consider the parallel primitivoids and the slant primitivoids of
$r\mbox{\rm -}\mathcal{P}r_{\bm {\gamma}}(t)$ and $\mathcal{P}r[\phi]_{\bm{\gamma}}$, respectively.
The normal vector fields of the above frontals are denoted by
$\bm{\nu}_{r\mbox{\rm -}\mathcal{P}r_{\bm{\gamma}}}(t)$ and
$\bm{\nu}_{\mathcal{P}r[\phi]_{\bm{\gamma}}}(t)$, respectively.
We have the following lemma.
\begin{Lem}
Let $(\bm{\gamma},\bm{\nu}):I\lon (\R^{2}\setminus\{\bm{0}\})\times S^{1}$ be a Legendrian curve
such that $\langle
\bm{\gamma}(t),\bm{\nu}(t)\rangle\not=0.$ For $\phi \not=\pi/2+n\pi$ {\rm (}$n\in \mathbb{Z}${\rm )}, we have
\[
\left\langle r\mbox{\rm -}\mathcal{P}r_{\bm{\gamma}}(t),\bm{\nu}_{r\mbox{\rm -}\mathcal{P}r_{\bm{\gamma}}}(t)\right\rangle\not= 0\ \mbox{and}\ 
\left\langle \mathcal{P}r[\phi]_{\bm{\gamma}}(t),\bm{\nu}_{\mathcal{P}r[\phi]_{\bm{\gamma}}}(t)\right\rangle\not= 0.
\]
\end{Lem}
\demo
By the definition of the primitive of $\bm{\gamma},$ the tangent vector of $\mathcal{P}r_{\bm{\gamma}}(t)$
at $s$ is $J\bm{\gamma}(t).$ So that the unit normal is 
\[
\bm{\nu}_{\mathcal{P}r_{\bm{\gamma}}}(t)=\pm JJ\frac{\bm{\gamma}}{\|\bm{\gamma}\|}(t)=
\pm \frac{\bm{\gamma}}{\|\bm{\gamma}\|}(t).
\]
It follows that
\[
\left\langle \mathcal{P}r_{\bm{\gamma}},\bm{\nu}_{\mathcal{P}r_{\bm{\gamma}}}\right\rangle
=\left\langle 2\bm{\gamma}-\frac{\|\bm{\gamma}\|^2}{\langle \bm{\nu},\bm{\gamma}\rangle}\bm{\nu},
\pm \frac{\bm{\gamma}}{\|\bm{\gamma}\|}\right\rangle=\pm \|\bm{\gamma}\|\not= 0.
\]
Since $r\mbox{\rm -}\mathcal{P}r_{\bm {\gamma}}(t)=r\mathcal{P}r_{\bm{\gamma}}(t),$
we have $\bm{\nu}_{r\mbox{\rm -}\mathcal{P}r_{\bm{\gamma}}}(t)=\pm \bm{\gamma}/\|\bm{\gamma}\|(t),$
so that the first assertion holds.
If we set $\bm{\nu}_{\phi}(t)=R(\phi)\bm{\nu}(t),$ then we have
$\langle \dot{\bm{\gamma}}_{\phi},\bm{\nu}_{\phi}\rangle =\langle R(\phi)\dot{\bm{\gamma}},R(\phi)\bm{\nu}\rangle=\langle \dot{\bm{\gamma}},\bm{\nu}\rangle=0.$
By definition, we can show that $\mathcal{P}r[\phi]_{\bm{\gamma}}=\cos\phi \mathcal{P}r_{\bm{\gamma}_\phi}.$
Since $\langle \bm{\gamma}_\phi,\bm{\nu}_\phi\rangle \not=0,$ we have
\[
\left\langle\mathcal{P}r[\phi]_{\bm{\gamma}},\bm{\nu}_{\mathcal{P}r[\phi]_{\bm{\gamma}}}\right\rangle=
\left\langle \cos\phi \mathcal{P}r_{\bm{\gamma}_\phi}, \bm{\nu}_{\mathcal{P}r_{\bm{\gamma}_\phi}}\right\rangle=\cos\phi \left\langle \mathcal{P}r_{\bm{\gamma}_\phi}, \bm{\nu}_{\mathcal{P}r_{\bm{\gamma}_\phi}}\right\rangle\not=0.
\]
This completes the proof.
\enD
Then we have the following theorem.
\begin{Th} Let $(\bm{\gamma},\bm{\nu}):I\lon (\R^{2}\setminus\{\bm{0}\})\times S^{1}$ be a Legendrian curve
such that $\langle
\bm{\gamma}(t),\bm{\nu}(t)\rangle\not=0.$ Suppose that $\phi\not=\pi/2+n\pi$ for $n\in \mathbb{Z}.$ Then
we have
\[
\cos (\psi+\phi)\mathcal{P}r[\psi]_{\mathcal{P}r[\phi]_{\bm{\gamma}}}(t)=\cos\psi\cos\phi\mathcal{P}r[\psi+\phi]_{\bm{\gamma}}(t).
\]
\end{Th}
\demo
By Theorem 5.3 and Lemma 5.4, $\mathcal{P}r[\phi]_{\bm{\gamma}}$ is a frontal.
Therefore, $\mathcal{P}r[\psi]_{\mathcal{P}r[\phi]_{\bm{\gamma}}}$ is well-defined.
By definition, we have
\[
\mathcal{P}r[\psi]_{\mathcal{P}r[\phi]_{\bm{\gamma}}}(t)=\cos\psi R(\psi)\mathcal{P}r[\phi]_{\bm{\gamma}}(t)=\cos\psi \mathcal{P}r_{R(\psi)\mathcal{P}r[\phi]_{\bm{\gamma}}}(t).
\]
Here, we have
\[
R(\psi)\mathcal{P}r[\phi]_{\bm{\gamma}}(t)=R(\psi)\cos\phi R(\phi)\mathcal{P}r_{\bm{\gamma}}(t)
=\cos\phi R(\psi+\phi)\mathcal{P}r_{\bm{\gamma}}(t)=
\cos\phi \mathcal{P}r_{\bm{\gamma}_{(\psi+\phi)}}(t).
\]
If $\psi+\phi \not= \pi/2+n\pi,$ then
\[
\cos\psi \mathcal{P}r_{R(\psi)\mathcal{P}r[\phi]_{\bm{\gamma}}}(t)=
\cos\psi \mathcal{P}r_{\cos\phi \mathcal{P}r_{\bm{\gamma}_{(\psi+\phi)}}}(t)
=\frac{\cos\psi\cos\phi}{\cos (\psi+\phi)}\mathcal{P}r[\psi+\phi]_{\bm{\gamma}}(t),
\]
so that we have
\[
\cos (\psi+\phi)\mathcal{P}r[\psi]_{\mathcal{P}r[\phi]_{\bm{\gamma}}}(t)=\cos\psi\cos\phi\mathcal{P}r[\psi+\phi]_{\bm{\gamma}}(t).
\]

If $\psi+\phi=\pi/2+n\pi,$ then
the both sides are zero.
\enD
\begin{Exa} We give an example of a front and its primitive.
We consider a curve defined by
\[
\bm{\gamma}(\theta)=\left(\frac{1}{32}(30 \cos\theta-17\cos3\theta+3\cos 5\theta),
\frac{1}{16\sqrt{2}}(\sin\theta (23+4\cos 2\theta -3\cos 4\theta))\right).
\]
We can show that this curve is a front, which is drawn in Fig. 7.
\begin{center}
\begin{tabular}{cc}   
\includegraphics[height=50mm,angle=90]{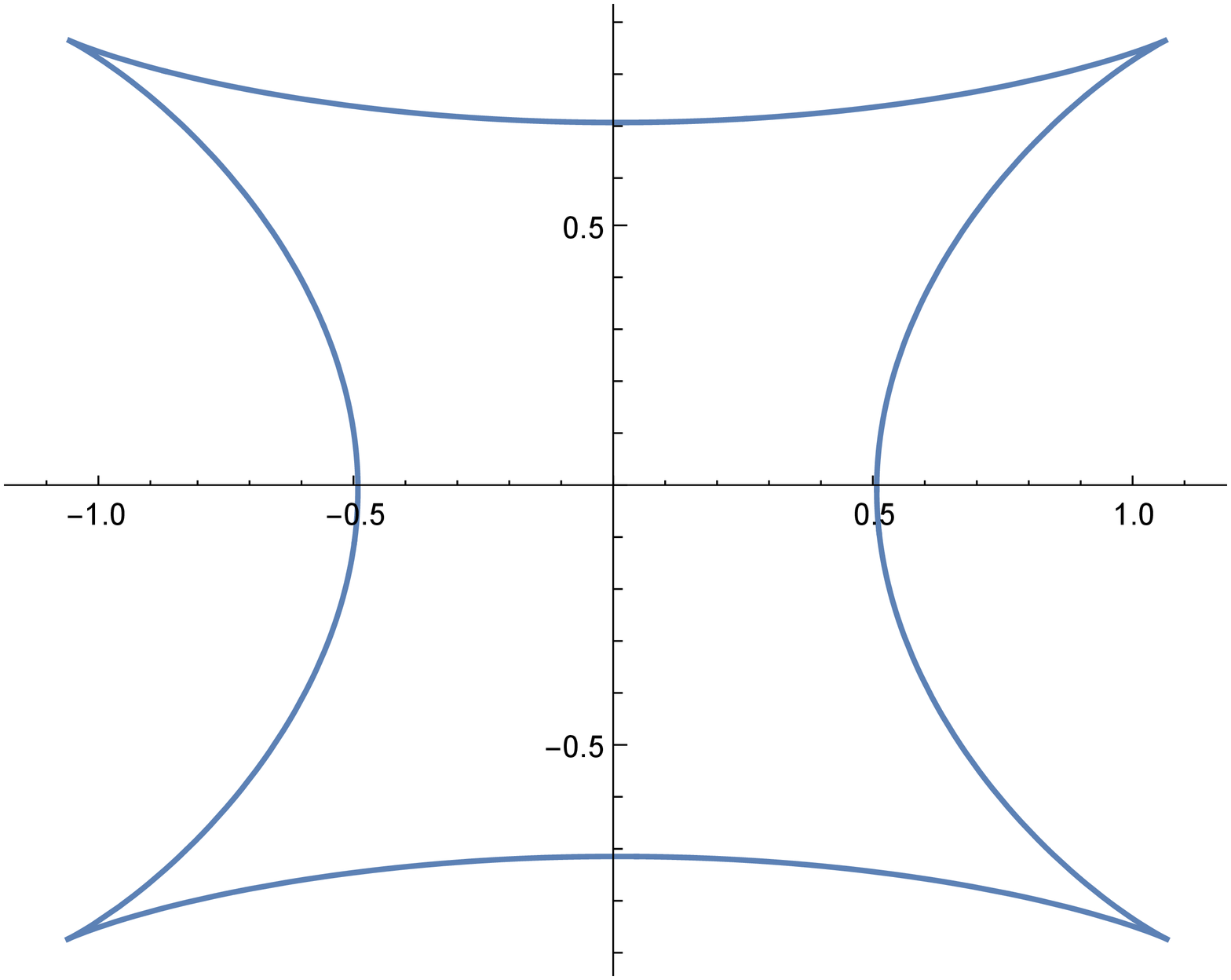}
&
\includegraphics[height=50mm,angle=90]{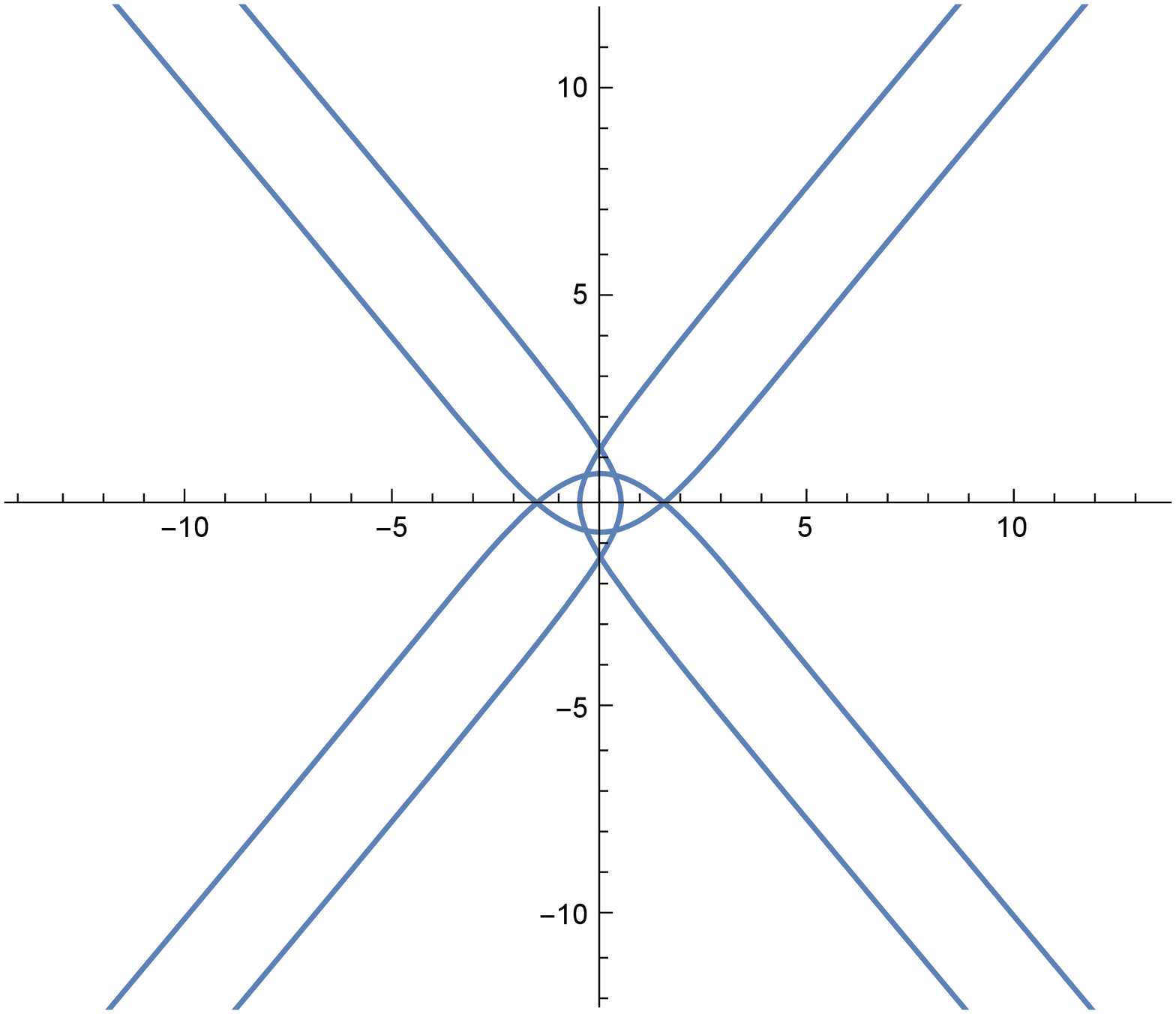}
 \\
\text{Fig. 7 : $\bm{\gamma}(\theta)$. }
&
\text{Fig. 8 : The primitive of $\bm{\gamma}$.}
\end{tabular}
\end{center}
The picture of the primitive of $\bm{\gamma}$ is given in Fig. 8.
The both of the pictures of $\bm{\gamma}$ and the primitive in Fig. 9.
We also draw the picture of the family of lines whose envelope is the primitive in Fig. 10.
of $\bm{\gamma}.$
\begin{center}
\begin{tabular}{cc}   
\includegraphics[height=50mm,angle=90]{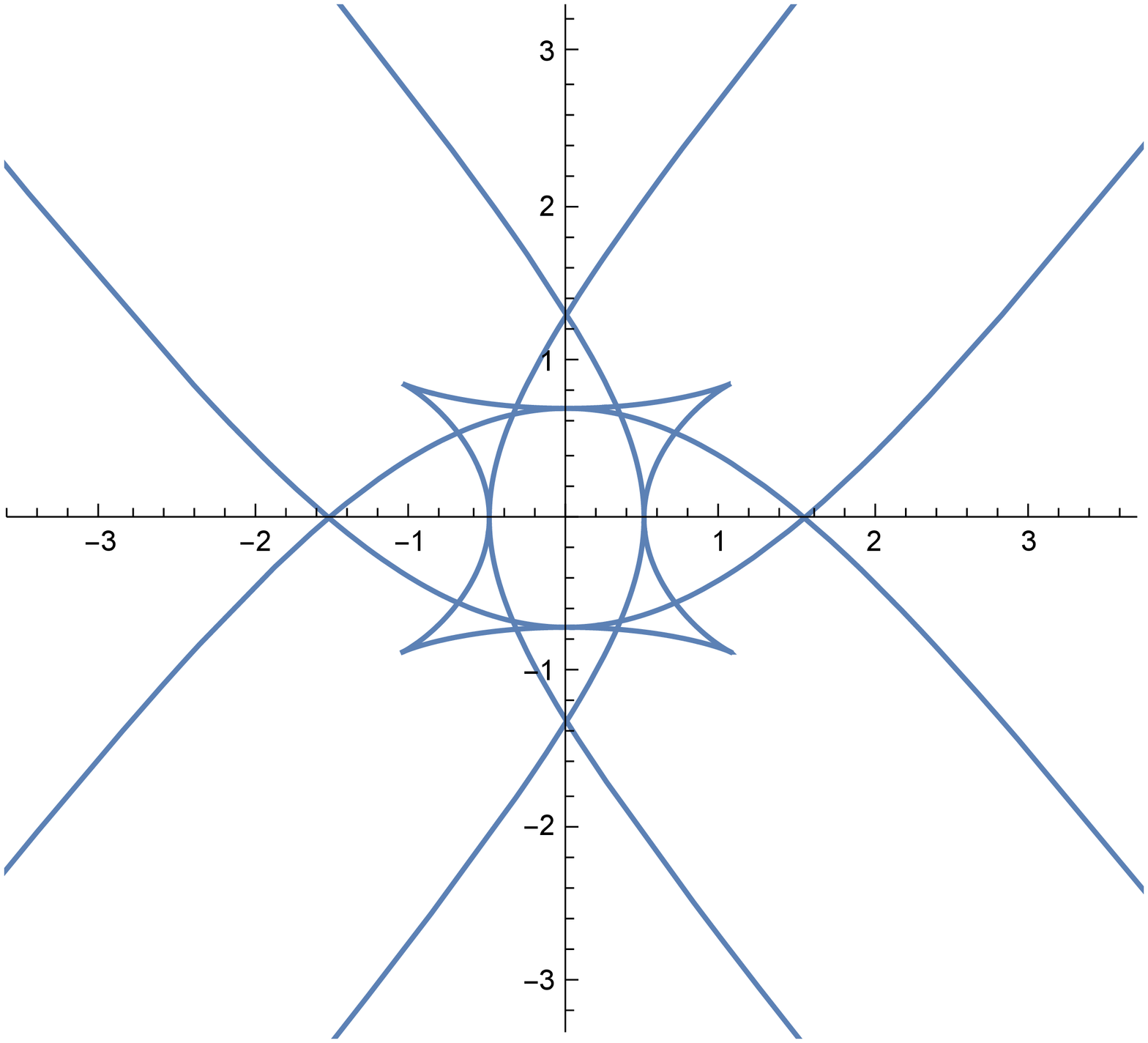}
&
\includegraphics[height=50mm,angle=90]{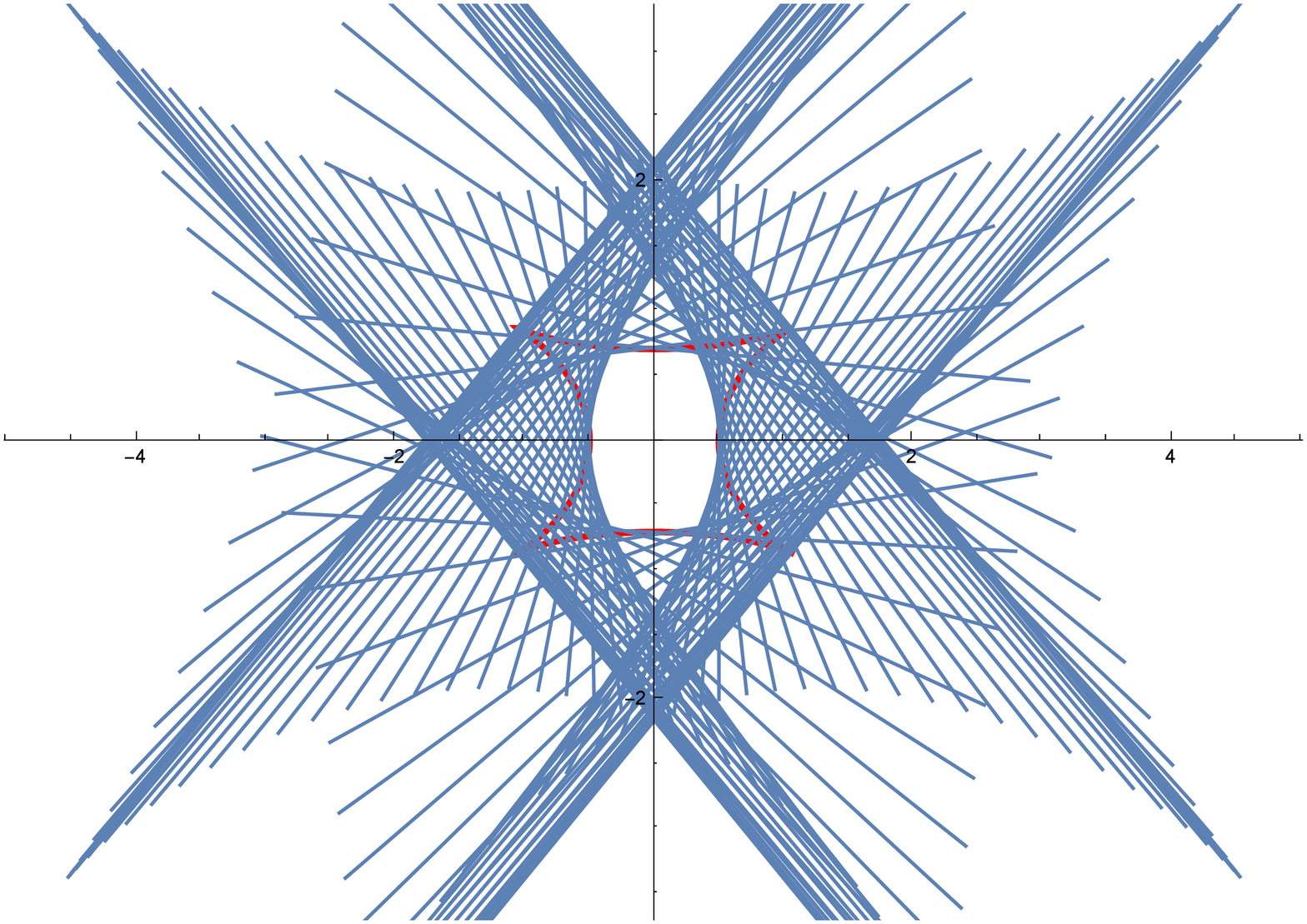}
 \\
\text{Fig. 9 : $\bm{\gamma}(\theta)$ and the primitive. }
&
\text{Fig. 10 : The primitive of $\bm{\gamma}$ with the family of lines.}
\end{tabular}
\end{center}
We can observe the aperture of the family of lines.
 Since pictures of all primitivoids are similar to the primitive, we omit to
 draw these pictures here.
We only remark that apertures of slant primitivoids always shrink to
the origin. This fact is related to the shape of distorted iris diaphragms of cameras.
We will discuss the relationship between iris diaphragms, evolutoids (cf. \cite{G-W}) and
slant primitivoids in elsewhere.

\end{Exa}

\begin{flushright}
\begin{tabular}{l}
Shyuichi Izumiya\\
Department of Mathematics, \\
Hokkaido University,\\
Sapporo 060-0810, Japan \hspace*{26mm}\\
{\tt izumiya@math.sci.hokudai.ac.jp}
\end{tabular}
\end{flushright}
\begin{flushright}
\begin{tabular}{l}
Nobuko Takeuchi\\
Department of Mathematics,\\
Tokyo Gakugei University,\\
Koganei, Tokyo, 184-8501, Japan \hspace*{12mm}\\
{\tt nobuko@u-gakugei.ac.jp}\\
\end{tabular}
\end{flushright}

\end{document}